\def\E{\end{document}}
\documentclass[11pt]{article}
\usepackage{amssymb,amsmath}

\begin{document}
\title{\bf
Global existence, blowup phenomena, and asymptotic behavior for quasilinear Schr\"{o}dinger equations
}
  \author{Xiaowei An$^{1,2}${\thanks{E-mail: anxiaowei2024@163.com (X. W. An).
Xiaowei An is supported by Key Research Priority Special Project of China People's Police University(ZDZX202601) and Hebei Provincial Natural Science Foundation(G2024507002).
}}\\
\small 1 School of Intelligence Policing, China People's Police University,\\
\small Langfang, He Bei, 065000, P. R. China\\
\small 2 Hebei Key Laboratory of Information Support Technology for Smart Policing,\\
\small China People's Police University, Langfang 065000, P. R. China\\
Xianfa Song$^{3,4}${\thanks{E-mail:\ \tt songxianfa@tju.edu.cn (or songxianfa2004@163.com)
  }}\\
\small 3 Department of Mathematics, School of Mathematics, Tianjin University,\\
\small Tianjin, 300072, P. R. China\\
\small 4 Xinjiang Production and Construction Corps\\
\small Key Laboratory of Green and Intelligent Development and Efficient Utilization of Strategic Mineral Resources,\\
\small Xinjiang University of Technology,  Hotan  Xinjiang 84800, P.R. China\\
}

\maketitle
\date{}

\newtheorem{theorem}{Theorem}[section]
\newtheorem{definition}{Definition}[section]
\newtheorem{lemma}{Lemma}[section]
\newtheorem{proposition}{Proposition}[section]
\newtheorem{corollary}{Corollary}[section]
\newtheorem{remark}{Remark}[section]
\renewcommand{\theequation}{\thesection.\arabic{equation}}
\catcode`@=11 \@addtoreset{equation}{section} \catcode`@=12

\begin{abstract}

In this paper, we study the Cauchy problem of the quasilinear Schr\"{o}dinger equation
\begin{equation*}
\left\{
\begin{array}{lll}
iu_t=\Delta u+2uh'(|u|^2)\Delta h(|u|^2)+F(|u|^2)u \quad {\rm for} \ x\in \mathbb{R}^N, \ t>0\\
u(x,0)=u_0(x),\quad x\in \mathbb{R}^N.
\end{array}\right.
\end{equation*}
Here $h(s)$ and $F(s)$ are some real-valued functions, with various choices for models from mathematical physics.
We examine the interplay between the quasilinear effect of $h$ and nonlinear effect of $F$ for the global existence and blowup phenomena.
We provide sufficient conditions on the blowup in finite time and global existence of the solution.
In some cases, we can deduce the watershed from these conditions. In the focusing case, we construct the sharp threshold for the blowup in finite time and global existence of the solution and lower bound for blowup rate of the blowup solution.



\end{abstract}

\section{Introduction}
\qquad In this paper, we consider the following Cauchy problem:
\begin{equation}
\label{1} \left\{
\begin{array}{lll}
iu_t=\Delta u+2uh'(|u|^2)\Delta h(|u|^2)+F(|u|^2)u \quad {\rm for} \ x\in \mathbb{R}^N, \ t>0\\
u(x,0)=u_0(x),\quad x\in \mathbb{R}^N.
\end{array}\right.
\end{equation}
Here $h(s)$ and $F(s)$ are some real valued functions, $h(s)\geq 0$ for $s\geq 0$. (\ref{1}) often appears in plasma physics and fluid mechanics, in the theory of Heisenberg ferromagnet and magnons, and in condensed matter theory, see \cite{Bass, Goldman, Litvak, Makhankov}. It has been used in many models of physical phenomena. For example, it models the self-channelling of a high-power ultra short laser in matter with $h(s)=\sqrt{1+s}$, the superfluid film  equation in plasma  physics  with $h(s)=s$ (\cite{Ku, LSS}) (so-called modified nonlinear Schr\"odinger equation),
and while with $h(s)=\sqrt{s}$, it illustrates the physical phenomenon in dissipative quantum mechanics, see \cite{Borovskii, Bouard, Hasse, Ritchie}.

The Cauchy problem of the equation $2i\frac{\partial a}{\partial t}+\Delta a+f(\gamma,\Delta\gamma)a=0$ subject to $a(x,0)=a_0(x)$  was considered by de Bouard, Hayashi and Saut in \cite{Bouard}. They obtained local as well as global existence results, where $\gamma=\sqrt{1+|a|^2}$. Later, the local well-posedness as well as the lower regularity of the Cauchy problem of quasilinear Schr\"odinger equation has been established by many authors, see  Theorem 1.1 in \cite{Colin1} as well as that in \cite{Colin2}, Theorem A and B in \cite{Kenig}, Theorem 1 in \cite{Marzuola1} as well as that in \cite{Marzuola2}, Theorem 1.2 and 1.3 in \cite{Marzuola3}, Theorem 6.1--Theorem 6.4 in \cite{Poppenberg1}, By these results, we have the following local existence result on (\ref{1}).

{\bf Theorem A}(Local well-posedness) {\it Assume that  $u_0\in H^{L+2}(\mathbb{R}^N)\cap X$ and $h(s)$, $F(s)\in C^{L+2}(\mathbb{R}^+)$ for some $L\geq N+2$. Then there exist a $T_L>0$ and a unique solution to (\ref{1}) satisfying
\begin{align*}
u(x,0)=u_0(x),\ u\in L^{\infty}([0,T_L];H^{L+2}(\mathbb{R}^N)\cap X)\cap C([0,T_L];H^L(\mathbb{R}^N)\cap X).
\end{align*}
Here
\begin{align}
X=\{w\in H^1(\mathbb{R}^N),\quad \int_{\mathbb{R}^N}|\nabla h(|w|^2)|^2dx<+\infty\}.\label{kongjianziji}
\end{align}
}

In order to state other results on (\ref{1}), we give the precise definition of global existence and finite time blowup of solutions.

{\bf Definition 1.} {\it Assume that $u(x,t)$ is the solution of (\ref{1}). If the maximum existence interval of $u(x,t)$ for $t$ is $[0, +\infty)$, we say that $u(x,t)$ is of global existence. On the other hand, we say that $u(x,t)$ blows up in finite time if there exists a time $0<T<+\infty$ such that
\begin{align}
\lim_{t\rightarrow T^-} \int_{\mathbb{R}^N}[|u(x,t)|^2+|\nabla u(x,t)|^2+|\nabla h(|u(x,t)|^2)|^2)]dx=+\infty.
\end{align}
 }


 One of main goals of this paper is to establish the global existence and blowup phenomena for the general problem (\ref{1}). About the topics on the global existence and blowup phenomena of the classical nonlinear Schr\"{o}dinger equation, in his seminal paper \cite{Glassey} of 1977, Glassey considered the following Cauchy problem
\begin{equation}
\label{2} \left\{
\begin{array}{lll}
iu_t=\Delta u+F(|u|^2)u \quad {\rm for} \ x\in \mathbb{R}^N, \ t>0\\
u(x,0)=u_0(x),\quad x\in \mathbb{R}^N.
\end{array}\right.
\end{equation}
By his result, the key condition on the blowup of the solution to (\ref{2}) is that there exists a constant $c_N>1+\frac{2}{N}$ such that $sF(s)\geq c_N G(s)$ for all $s\geq 0$, where $G(s)=\int_0^s F(\eta)d\eta$. In 1981, Berestycki and Cazenave,  Cazenave and Lions established a sharp condition on the blowup of the solution and the result on the instability of standing wave solution to (\ref{2}) in \cite{Berestycki, Cazenave}. We also refer to \cite{Coz, Zhang1} and the references therein.  In \cite{Guo}, Guo, Chen and Su, studied the following problem:
\begin{equation}
\label{xj} \left\{
\begin{array}{lll}
i\varphi_t+\Delta \varphi+2(\Delta |\varphi|^2)\varphi+|\varphi|^{q-2}\varphi=0\quad {\rm for} \ x\in \mathbb{R}^N, \ t>0\\
\varphi(x,0)=\bar{u}_0(x),\quad x\in \mathbb{R}^N.
\end{array}\right.
\end{equation}
Letting $\bar{\varphi}=u$, we find that (\ref{xj}) is a special case of (\ref{1}) with $h(s)=s$ and $F(s)=s^{\frac{q-2}{2}}$. They obtained that
the solution of (\ref{xj}) will blow up in finite time if $4+\frac{4}{N}\leq q<2\cdot 2^*$ under some other assumptions. Here $2^*=\frac{2N}{N-2}$.

About the existence of standing wave solution to (\ref{1}), we can refer to \cite{Chen1, Colin2, Deng, Fang, Jeanjean, Liu1, Liu2, Liu3, Poppenberg2, Shen} and the references therein. Recently, the stability and instability of the standing wave solution of (\ref{1}) with $h(s)=s^{\alpha}$ and $F(s)=s^{\frac{q-1}{2}}$ was respectively studied by Colin, Jeanjean and Squassina in \cite{Colin2}(where $\alpha=1$),  Chen, Li and Wang in \cite{Chen1}(where $\alpha>\frac{1}{2}$).
Letting $\varphi=\bar{u}$, the models in \cite{Colin2} and \cite{Chen1} can be written as
\begin{equation}
\label{3} \left\{
\begin{array}{lll}
i\varphi_t+\Delta \varphi+2\alpha(\Delta |\varphi|^{2\alpha})|\varphi|^{2\alpha-2}\varphi+|\varphi|^{q-2}\varphi=0\quad {\rm for} \ x\in \mathbb{R}^N, \ t>0\\
\varphi(x,0)=\bar{u}_0(x),\quad x\in \mathbb{R}^N.
\end{array}\right.
\end{equation}
By their results, the standing wave solution of (\ref{3}) is stable if $2<q<4\alpha+\frac{4}{N}$ and unstable if $4\alpha+\frac{4}{N}\leq q<2\alpha\cdot 2^*$. Chen and Rocha in \cite{Chen2} studied the equation with a harmonic potential
\begin{equation}
\label{3xj} \left\{
\begin{array}{lll}
i\varphi_t+\Delta \varphi+2(\Delta |\varphi|^2)\varphi-|x|^2\varphi+|\varphi|^{q-2}\varphi=0\quad {\rm for} \ x\in \mathbb{R}^N, \ t>0\\
\varphi(x,0)=\bar{u}_0(x),\quad x\in \mathbb{R}^N.
\end{array}\right.
\end{equation}
They proved the standing wave solution of (\ref{3xj}) is stable if $2<q<4+\frac{4}{N}$ and unstable if $4+\frac{4}{N}\leq q<2\cdot 2^*$.

Motivated by these studies, we investigate problem (\ref{1}) for general quasilinear term and nonlinearity. We will establish conditions on the blowup in finite time and global existence of the solution to the more general equation (\ref{1}). Our results reveal the crucial interplay between $h$ and $F$ for the behavior of solutions structure.

Before we state our results, we define the mass and energy of (\ref{1}) as follows.

(i) Mass: $$ m(u)=\left(\int_{\mathbb{R}^N}|u(\cdot,t)|^2dx\right)^{\frac{1}{2}}:=[M(u)]^{\frac{1}{2}};$$

(ii) Energy : $$E(u)=\frac{1}{2}\int_{\mathbb{R}^N}[|\nabla u|^2+|\nabla h(|u|^2)|^2-G(|u|^2)]dx.$$

We will prove the conservations of mass and energy in Section 2.

In the sequels, we will use $C$, $C'$, $c_j$, $c'_j$, and so on, to denote the constants which are independent of $x$ and $t$, the values of them may vary line to line.

We use $C_s$ to denote the best constant in the Sobolev's inequality
\begin{align}
\int_{\mathbb{R}^N}w^{2^*}dx\leq C_s\left(\int_{\mathbb{R}^N}|\nabla w|^2dx\right)^{\frac{2^*}{2}}\quad {\rm for \ any}\quad w\in H^1(\mathbb{R^N}),\label{zjcs}
\end{align}

Our first result establishes sufficient conditions on the global existence of the solution to (\ref{1}).

{\bf Theorem 1.} {\it  Let $u(x,t)$ be the solution to (\ref{1}) with $u_0\in X$. Assume $F(s)=F_1(s)-F_2(s)$,
 and denote
$$G_1(s)=\int_0^sF_1(\eta)d\eta, \quad G_2(s)=\int_0^sF_2(\eta)d\eta.$$

Then $\int_{\mathbb{R}^N}[|u|^2+|\nabla u|^2+|\nabla h(|u|^2)|^2]dx$ is uniformly bounded for all $t>0$(i.e., $u$ is of global existence) in one of the following three cases:

Case (1) $F_1(s)\equiv 0$, $F(s)\equiv -F_2(s)\leq 0$ for $s\geq 0$;

Case (2) $F_2(s)\equiv 0$, $F(s)\equiv F_1(s)$. Suppose that $F_1(s)\geq 0$ for $s\geq 0$ or $F_1(s)$ changes sign for $s\geq 0$, and there exist $0<\theta_1<1$, $0<\theta_2<1$, $q_1>1$ and $q_2>1$ such that
\begin{align}
(2^*-2)\theta_1+2q_1\geq 2^*,\quad (2^*-2)\theta_2+2q_2\geq 2^*,\label{1015x2}
\end{align}
\begin{align}
&[|G_1(s)|]^{\theta_1}\leq c_1s,\quad [|G_1(s)|]^{q_1}\leq c'_1s+\epsilon_1[s^{\frac{1}{2}}+h(s)]^{2^*}\quad {\rm for} \quad\ 0\leq s\leq 1, \label{1015x1'}\\
&  [|G_1(s)|]^{\theta_2}\leq c_2s,\quad
[|G_1(s)|]^{q_2}\leq c'_2s+\epsilon_2[s^{\frac{1}{2}}+h(s)]^{2^*}\quad {\rm for} \quad s>1 \label{1015x1}
\end{align}
 for some positive constants $c_1, c'_1, c_2, c'_2, \epsilon_1$. Moreover, the initial value $u_0$ satisfies

(i) $$2^{\frac{2^*-1}{\tau'_1}}c_1^{\frac{1}{\tau_1}}(\epsilon_1C_s)^{\frac{1}{\tau'_1}}\|u_0\|_2^{\frac{2}{\tau_1}}<1$$
if
\begin{align}
(2^*-2)\theta_1+2q_1=2^*,\quad (2^*-2)\theta_2+2q_2>2^*,\label{615xj1}
\end{align}

(ii) $$2^{\frac{2^*-1}{\tau'_2}}c_2^{\frac{1}{\tau_2}}(\epsilon_2C_s)^{\frac{1}{\tau'_2}}\|u_0\|_2^{\frac{2}{\tau_2}}<1$$
if
\begin{align}
(2^*-2)\theta_1+2q_1>2^*,\quad (2^*-2)\theta_2+2q_2= 2^*, \label{615xj2}
\end{align}

(iii) $$2^{\frac{2^*-1}{\tau'_1}}c_1^{\frac{1}{\tau_1}}(\epsilon_1C_s)^{\frac{1}{\tau'_1}}\|u_0\|_2^{\frac{2}{\tau_1}}
+2^{\frac{2^*-1}{\tau'_2}}c_2^{\frac{1}{\tau_2}}(\epsilon_2C_s)^{\frac{1}{\tau'_2}}\|u_0\|_2^{\frac{2}{\tau_2}}<1$$ if
\begin{align}
(2^*-2)\theta_1+2q_1=2^*,\quad (2^*-2)\theta_2+2q_2=2^*. \label{615xj1'}
\end{align}
Here
$$
\frac{1}{\tau_j}=\frac{q_j-1}{q_j-\theta_j},\qquad \frac{1}{\tau'_j}=\frac{1-\theta_j}{q_j-\theta_j},\quad j=1,2;
$$

Case (3) $F_1(s)\geq 0, \not\equiv 0$ or $F_1(s)$ changes sign for $s\geq 0$, $F_2(s)\geq 0, \not\equiv 0$ for $s\geq 0$.
Suppose that

(iv) $F_1(s)$ satisfies the assumptions of Case (2) or

(v) there exists $c$ such that
\begin{align}
|G_1(s)|\leq cs+G_2(s)\quad {\rm for} \ s\geq 0\label{6241}
\end{align}
or

(vi) there  exist $0<\tilde{\alpha}_1<1$, $0<\tilde{\alpha}_2<1$,  $\tilde{\beta}_1>1$ and $\tilde{\beta}_2>1$ such that
\begin{align}
&[|G_1(s)|]^{\tilde{\alpha}_1}\leq \tilde{c}_1s,\quad [|G_1(s)|]^{\tilde{\beta}_1}\leq \tilde{c}_2G_2(s)\quad {\rm for}\ 0\leq s\leq 1,\label{6171}\\
&[|G_1(s)|]^{\tilde{\alpha}_2}\leq \tilde{c}'_1s,\quad [|G_1(s)|]^{\tilde{\beta}_2}\leq \tilde{c}'_2G_2(s)\quad {\rm for}\  s>1,\label{6172}
\end{align}
for some constants $\tilde{c}_1>0$, $\tilde{c}'_1>0$, $\tilde{c}_2>0$ and $\tilde{c}'_2>0$.

}

{\bf Remark 1.1.} 1. If $h(s)\equiv 0$, $F(s)=bs^{\tilde{q}}(b>0)$, we can take $\theta_1=\theta_2=\frac{1}{\tilde{q}+1}$, $q_1=\frac{(\tilde{q}+1)\cdot 2^*+1}{\tilde{q}+1}$, $q_2=\frac{2^*}{2(\tilde{q}+1)}$, and the solution is global existence when $\tilde{q}<\frac{2}{N}$.
Our results meet with the classic ones on semilinear Schr\"{o}dinger equation in \cite{Ginibre3, Ginibre4}.

2. If $h(s)=as^{p}(a>0)$, $F(s)=bs^{\tilde{q}}(b>0)$. We can take $\theta_1=\theta_2=\frac{1}{\tilde{q}+1}$, $q_1=\frac{(\tilde{q}+1)\cdot 2^*+1}{\tilde{q}+1}$. $q_2=\frac{p\cdot 2^*}{\tilde{q}+1}$ if $p\geq \frac{1}{2}$ and $q_2=\frac{\cdot 2^*}{2(\tilde{q}+1)}$ if $p\leq \frac{1}{2}$. Meanwhile, the conditions $q_2>1$ and $(2^*-2)\theta_2+2q_2>2^*$ imply that $0<\tilde{q}<\max\{\frac{2}{N}, 2p-1+\frac{2}{N}\}$.
Our results meet with those on quasilinear Schr\"{o}dinger equation in \cite{Chen1, Chen2, Colin2, Guo}.\hfill $\Box$

Our second result will establish the sufficient conditions on the blowup in finite time for the solutions to (\ref{1}).

{\bf Theorem 2.} {\it Let $u(x,t)$ be the solution to (\ref{1}) with $u_0\in X$. Assume that there exist constants $k$, $c_N$ and
$c_M$ such that

(i) $sh''(s)\leq kh'(s)$ if $h'(s)\geq 0$ or  $sh''(s)\geq kh'(s)$ if $h'(s)\leq 0$;

(ii) $c_M\geq 0$, $c_N>\max\{1+\frac{2}{N}, 2(k+1)+\frac{2}{N}\}$ and $c_NG(s)\leq sF(s)+c_Ms$.

Suppose that $\Im \int_{\mathbb{R}^N}\bar{u}_0(x\cdot \nabla u_0)dx>0$, $|x|u_0\in L^2(\mathbb{R}^N)$, $$2(c_N-1)E(u_0)+c_MM(u_0)\leq 0\qquad {\rm if}\quad k\leq -\frac{1}{2}$$  or $$2[(2k+1)N+2]E(u_0)+c_MM(u_0)\leq 0 \qquad {\rm if} \quad k>-\frac{1}{2}.$$

Then there exists a finite time $T$ such that
$$
\lim_{t\rightarrow T^-} \int_{\mathbb{R}^N}[|u(x,t)|^2+|\nabla u(x,t)|^2+|\nabla h(|u|^2)(x,t)|^2]dx=+\infty.
$$
}

As direct consequences of Theorem 2, we give two corollaries, which answer the question of how to determine the constants $k$, $c_N$ and $c_M$ in relation to $h(s)$ and $F(s)$.

{\bf Corollary 1.1. }{\it Let $u(x,t)$ be the solution to (\ref{1}) with $u_0\in X$, and  $F(s)\geq 0$ for $s\geq 0$ and
$$
k= \inf_{\tilde{k}}\{h'(s)\geq 0,\quad sh''(s)\leq \tilde{k}h'(s) \quad  {\rm for}\ s\geq 0 \}.
$$
Suppose that

(i) $E(u_0)\leq 0$;

(ii) $\Im \int_{\mathbb{R}^N}\bar{u}_0(x\cdot \nabla u_0)dx>0$, $|x|u_0\in L^2(\mathbb{R}^N)$;

(iii) There exists  $c_N>\max\{1+\frac{2}{N}, 2(k+1)+\frac{2}{N}\}$ such that $c_NG(s)\leq sF(s)$ for $s\geq 0$.

Then there exists a finite time $T$ such that
$$
\lim_{t\rightarrow T^-} \int_{\mathbb{R}^N}[|u(x,t)|^2+|\nabla u(x,t)|^2+|\nabla h(|u|^2)(x,t)|^2]dx=+\infty.
$$
}

{\bf Corollary 1.2. }{\it Let $u(x,t)$ be the solution to (\ref{1}) with $u_0\in X$, and  $F(s)\geq 0$ for $s\geq 0$ and
$$
k= \inf_{\tilde{k}}\{h'(s)\geq 0,\quad sh''(s)\leq \tilde{k}h'(s) \quad  {\rm for}\ s\geq 0 \}.
$$
Assume that there exist $c_N$ and $0\leq c_M<\bar{c}_M$ such that
\begin{align*}
 c_N=\sup\{\bar{c}_N: \bar{c}_NG(s)\leq sF(s)+c_Ms\ {\rm for}\ s\geq 0\}>\bar{C}(k,N)\}.
\end{align*}
Here
$$
\bar{C}(k,N)= \max\{1+\frac{2}{N},\ 2(k+1)+\frac{2}{N}\}
$$
and
$$
\bar{c}_M=\inf\{\tilde{c}_M: \max\{2^*,\ 2k+1+2^*\}G(s)\leq sF(s)+\tilde{c}_Ms\}.
$$
Suppose that $\Im \int_{\mathbb{R}^N}\bar{u}_0(x\cdot \nabla u_0)dx>0$, $|x|u_0\in L^2(\mathbb{R}^N)$, $2(c_N-1)E(u_0)+c_MM(u_0)\leq 0$ if $k\leq -\frac{1}{2}$ or $2[(2k+1)N+2]E(u_0)+c_MM(u_0)\leq 0$ if $k>-\frac{1}{2}$.

Then there exists a finite time $T$ such that
$$
\lim_{t\rightarrow T^-} \int_{\mathbb{R}^N}[|u(x,t)|^2+|\nabla u(x,t)|^2+|\nabla h(|u|^2)(x,t)|^2]dx=+\infty.
$$
}

{\bf Remark 1.2.} First, we would like to compare our results on (\ref{1}) with those on semilinear Sch\"{o}dinger equation. If $h(s)\equiv 0$, the first equation in (\ref{1}) is a semilinear Sch\"{o}dinger equation, we can take $k=-\frac{1}{2}$ and $c_M=0$ in Theorem 2, our result meets with the classic results on semilinear Schr\"{o}dinger equation in \cite{Cazenave2003, Glassey}.

If $c_N>\max\{1+\frac{2}{N}, 2(k+1)+\frac{2}{N}\}$, then the solution of (\ref{1}) blows up in finite time under $E(u_0)\leq 0$ and other conditions. Naturally, if $c_N>\max\{1+\frac{2}{N}, 2(k+1)+\frac{2}{N}\}$ and $E(u_0)>0$, it is interesting to know whether the solution is of global existence or of blowing up in finite time. The third result  answers this question and establishes a sharp threshold on the blowup and global existence of the solution to (\ref{1}).

{\bf Remark 1.3.}  We also give some examples of $h(s)$ and $F(s)$ to illustrate the results of Theorem 2 below.

1. $h(s)=as^{p}(a>0)$, $F(s)=bs^{\tilde{q}}(b>0)$. If $0<p$, then we can take $k=(p-1)$, $c_N=\tilde{q}+1$, $c_M=0$ and under the assumptions (i) and (ii), the solution will blow up in finite time when $\tilde{q}=c_N-1>\max\{\frac{2}{N}, 2p-1+\frac{2}{N}\}$ by Theorem 2.  Especially, if $p=1$, then we can take $k=0$, and the solution will blow up in finite time when $\tilde{q}>1+\frac{2}{N}$, our result meets with that of \cite{Guo}. On the other hand, the solution is global existence for any initial data when $\tilde{q}=c_N-1<\max\{\frac{2}{N}, 2p-1+\frac{2}{N}\}$ by Theorem 1.

2. $h(s)=a_1s^{p_1}+a_2s^{p_2}+...+a_ms^{p_m}$, $F(s)=b_1s^{\tilde{q}_1}+b_2s^{\tilde{q}_2}+...+b_ns^{\tilde{q}_n}$, $1<p_1<p_2<...<p_m$, $a_l>0$, $l=1,2,...,m$, $0<\tilde{q}_1<\tilde{q}_2<...<\tilde{q}_n$, $b_j>0$,  $j=1,2,...,n$.  We can take $k=p_m-1$ and $c_N=\tilde{q}_n+1-\epsilon$, $c_M=C(\tilde{q}_1,...,\tilde{q}_n, b_1,...,b_n,\epsilon)$ for $0<\epsilon<\tilde{q}_n-(2p_m-1+\frac{2}{N})$, and under the assumptions (i) and (ii), the solution will blow up in finite time if $\tilde{q}_n>2p_m-1+\frac{2}{N}$ by Theorem 2. On the other hand, and the solution is global existence if $0<\tilde{q}_n<2p_m-1+\frac{2}{N}$ by Theorem 1.

3. $h(s)=as^{p}(a>0)$, $F(s)=b_1s^{\tilde{q}_1}-b_2s^{\tilde{q}_2}$, $b_1>0$, $b_2>0$, $\tilde{q}_1>\tilde{q}_2>0$, $G(s)=\frac{b_1}{\tilde{q}_1+1}s^{\tilde{q}_1+1}-\frac{b_2}{\tilde{q}_2+1}s^{\tilde{q}_2+1}$, $(\tilde{q}_1+1)G(s)\leq sF(s)$, then we can take $k=(p-1)$, $c_N=\tilde{q}+1$, $c_M=0$ and under the assumptions (i) and (ii), the solution will blow up in finite time when $\tilde{q}_1=c_N-1>\max\{\frac{2}{N}, 2p-1+\frac{2}{N}\}$ by Theorem 2. On the other hand, if $\max\{\frac{2}{N}, 2p-1+\frac{2}{N}\}>\tilde{q}_1>\tilde{q}_2>0$ or $0<\tilde{q}_1<\tilde{q}_2$,
the solution is global existence by Theorem 1.

{\bf Theorem 3.} (Sharp Threshold ) {\it Let $u(x,t)$ be the solution of (\ref{1}) with $u_0\in X$ and $F(s)\geq 0$ for $s\geq 0$.

 Assume that (i) there exist constants $k$, $c_N$, $2\leq L<N(c_N-1)$ and $0<\underline{l}\leq 2$ such that
\begin{align*}
& sh''(s)\leq kh'(s)\quad  {\rm if} \ h'(s)\geq 0\quad {\rm or} \  sh''(s)\geq kh'(s) \quad {\rm if} \  h'(s)\leq 0,\\
& (L-2)+4[L-(N+2)][h'(s)]^2s-8Nh''(s)h'(s)s^2\geq 0,\\
&(2-\underline{l})+4(N+2-\underline{l})[h'(s)]^2s+8Nh''(s)h'(s)s^2\geq 0
\end{align*}
and  $\quad c_NG(s)\leq sF(s)$ for $s\geq 0$;

(ii) There exist $0<\theta_3<1$, $0<\theta_4<1$, $q_3>1$ and $q_4>1$ such that
\begin{align}
(2^*-2)\theta_3+2q_3\leq 2^*,\quad (2^*-2)\theta_4+2q_4\leq 2^*\label{1017s2}
\end{align}
and
\begin{align}
&[sF(s)]^{\theta_3}\leq c_3s,\quad [sF(s)]^{q_3}\leq c'_3[s^{\frac{1}{2}}+h(s)]^{2^*}\quad {\rm for} \quad\ 0\leq s\leq 1, \label{1017s1'}\\
&  [sF(s)]^{\theta_4}\leq c_4s,\quad
[sF(s)]^{q_4}\leq c'_4[s^{\frac{1}{2}}+h(s)]^{2^*}\quad {\rm for} \quad s>1. \label{1017s1}
\end{align}
with some positive constants $c_3, c'_3, c_4$ and $c'_4$.

Moreover, suppose that
there exists $\omega>0$ such that
\begin{align} d_I:=\inf_{\{w\in H^1(\mathbb{R}^N)\setminus \{0\};
Q(w)=0\}}\left(\frac{\omega}{2}\|w\|_2^2+E(w)\right)>0,\label{9651}\end{align}
where \begin{align}
Q(w)&=2\int_{\mathbb{R}^N}|\nabla w|^2dx+(N+2)\int_{\mathbb{R}^N}|\nabla h(|w|^2)|^2dx\nonumber\\
&\quad +8N\int_{\mathbb{R}^N}h''(|w|^2)h'(|w|^2)|w|^4|\nabla w|^2dx-N\int_{\mathbb{R}^N}[|w|^2F(|w|^2)-G(|w|^2)]dx,\\
E(w)&=\frac{1}{2}\int_{\mathbb{R}^N}[|\nabla w|^2+|\nabla h(|w|^2)|^2-G(|w|^2)]dx,
\end{align}
and $u_0$
satisfies
$$\frac{\omega}{2}\|u_0\|_2^2+E(u_0)<d_I.$$

Then we have:

(1). If $Q(u_0)>0$, the solution of (\ref{1}) exists
globally;

(2). If $Q(u_0)<0$ and $\Im \int_{\mathbb{R}^N} \bar{u}_0(x\cdot \nabla u_0)dx\geq 0$, $|x|u_0\in L^2(\mathbb{R}^N)$, the solution
of (\ref{1}) blows up in finite time.}

{\bf Remark 1.4.} 1. Under the assumptions of Theorem 3, by the results of \cite{Colin1,Deng,Shen}, the minimizer of (\ref{9651}) can be achieved at $w(x)$
which is a weak solution of
$$
-\Delta w-2wh'(w^2)\Delta h(w^2)+\omega w-F(w^2)w=0
$$ for any $\omega>0$. Hence $$Q(w)=0,\quad \omega \|w\|_2^2+E(w)=d_I.$$

2. Some examples of $h(s)$ and $F(s)$ which satisfy the assumptions of Theorem 3 will be given in Remark 5.2.

Having dealt with the conditions on  blowup in finite time and the global existence of the solutions to (\ref{1}), we will consider asymptotic behavior for
the solutions. Inspired by \cite{Cazenave2003, Ginibre1, Ginibre2}, we have the pseudo-conformal conservation law below, which is essential for the study of the asymptotic behavior of the global solutions and the lower bound for the  blowup rate of the blowup solution.
Let $u$ be a solution of (\ref{1}). We set
\begin{align}
\theta(t)&=\int_{\mathbb{R}^N}-4N[2h''(|u|^2)h'(|u|^2)|u|^2+( h'(|u|^2))^2]|u|^2|\nabla u|^2dx\nonumber\\
&\qquad+\int_{\mathbb{R}^N} [NF(|u|^2)|u|^2-(N+2)G(|u|^2)]dx.\label{691'}
\end{align}
{\bf Theorem 4.( Pseudo-conformal Conservation Law)}

 {\it 1. Assume that $u$ is the global solution of (\ref{1}), $u_0\in X$ and $xu_0\in L^2(\mathbb{R}^N)$. Let $G(s)=\int_0^sF(\eta)d\eta$. Then
\begin{align}
P(t)&=\int_{\mathbb{R}^N}|(x-2it\nabla)u|^2dx+4t^2\int_{\mathbb{R}^N}|\nabla h(|u|^2)|^2dx-4t^2\int_{\mathbb{R}^N}G(|u|^2)dx\nonumber\\
&=\int_{\mathbb{R}^N}|xu_0|^2dx+4\int_0^t\tau\theta(\tau)d\tau.\label{691}
\end{align}

2.  Assume that $u$ is the blowup solution of (\ref{1}) with blowup time $T$, $u_0\in X$ and $xu_0\in L^2(\mathbb{R}^N)$. Then
\begin{align}
B(t)&:=\int_{\mathbb{R}^N}|(x+2i(T-t)\nabla)u|^2dx+4(T-t)^2\int_{\mathbb{R}^N}|\nabla h(|u|^2)|^2dx-4(T-t)^2\int_{\mathbb{R}^N}G(|u|^2)dx\nonumber\\
&=\int_{\mathbb{R}^N}|(x+2iT\nabla)u_0|^2dx+4T^2\int_{\mathbb{R}^N}[|\nabla h(|u_0|^2)|^2-G(|u_0|^2)]dx\nonumber\\
&\quad-4\int_0^t(T-\tau)\theta(\tau)d\tau.\label{893}
\end{align}
}

As the applications of Theorem 4, in Section 6 we will give some asymptotic behavior results on the global solution of (\ref{1}) and the lower bound for the  blowup rate the blowup solution of  (\ref{1}) (see Theorem 5).

The organization of this paper is as follows. In Section 2, we will prove some equalities which will be applied to prove other conclusions later.
In Section 3, we will prove Theorem 1, which will establish the sufficient conditions on the global existence of the solution to (\ref{1}). In Section 4, we will prove Theorem 2, which will establish the sufficient conditions on  blowup in finite time for the solution to (\ref{1}). In Section 5, we will prove Theorem 3, which will establish a sharp threshold on the blowup in finite time and global existence of the solution to (\ref{1}). In Section 6, we will prove Theorem 4 and Theorem 5, which will give some asymptotic behavior results on the global solution of  (\ref{1}).

\section{Preliminaries}
\qquad In this section, we will prove a lemma as follows.

{\bf Lemma 2.1.} {\it Assume that $u$ is the solution to (\ref{1}). Then in the time interval $[0,t]$ when it exists, $u$ satisfies

(i) Mass conversation: $$ M(u)=\|u(\cdot,t)\|_2^2=M(u_0)=\|u_0\|^2_2;$$

(ii) Energy conversation: $$E(u)=\frac{1}{2}\int_{\mathbb{R}^N}[|\nabla u|^2+|\nabla h(|u|^2)|^2-G(|u|^2)]dx=E(u_0);$$

(iii) $$\frac{d}{dt} \int_{\mathbb{R}^N}|x|^2|u|^2dx=-4\Im \int_{\mathbb{R}^N} \bar{u}(x\cdot \nabla u)dx;$$

(iv) \begin{align*}
&\quad \frac{d}{dt} \Im \int_{\mathbb{R}^N} \bar{u}(x\cdot \nabla u)dx=-2\int_{\mathbb{R}^N}|\nabla u|^2dx-(N+2)\int_{\mathbb{R}^N}|\nabla h(|u|^2)|^2dx\\
&\qquad -8N\int_{\mathbb{R}^N}h''(|u|^2)h'(|u|^2)|u|^4|\nabla u|^2dx+N\int_{\mathbb{R}^N}[|u|^2F(|u|^2)-G(|u|^2)]dx,
\end{align*}
where $G(s)=\int_0^s F(\eta)d\eta$.

}

{\bf Proof:} (i) Multiplying (\ref{1}) by $2\bar{u}$, taking the imaginary part of the result, we have
\begin{align}
\frac{\partial }{\partial t}|u|^2=\Im(2\bar{u}\Delta u) =\nabla \cdot (2\Im \bar{u}\nabla u).\label{10121}
\end{align}
Integrating it over $\mathbb{R}^N\times [0,t]$, we have
$$ \int_{\mathbb{R}^N}|u|^2dx=\int_{\mathbb{R}^N}|u_0|^2dx.$$

(ii)  Multiplying (\ref{1}) by $2\bar{u}_t$, taking the real part of the result, then integrating it over $\mathbb{R}^N\times [0,t]$, we have
$$\int_{\mathbb{R}^N}[|\nabla u|^2+|\nabla h(|u|^2)|^2-G(|u|^2)]dx=\int_{\mathbb{R}^N}[|\nabla u_0|^2+|\nabla h(|u_0|^2)|^2-G(|u_0|^2)]dx.$$

(iii) Multiplying (\ref{10121}) by $|x|^2$ and integrating it over $\mathbb{R}^N$, we have
\begin{align*}
\frac{d}{dt}\int_{\mathbb{R}^N}|x|^2|u|^2dx&=\int_{\mathbb{R}^N}|x|^2\nabla \cdot(2\Im (\bar{u}\nabla u))dx
=-4\Im \int_{\mathbb{R}^N}\bar{u}(x\cdot \nabla u)dx.
\end{align*}

(iv) Denote $u(x,t)=a(x,t)+ib(x,t)$, i.e.,  $a(x,t)=\Re u(x,t)$ and $b(x,t)=\Im u(x,t)$. Then
$$
\frac{d}{dt}\Im \bar{u}(x\cdot \nabla u)=\sum_{k=1}^N[x_k(b_t)_{x_k}a-x_k(a_t)_{x_k}b]+\sum_{k=1}^N(x_kb_{x_k}a_t-x_ka_{x_k}b_t).
$$
And
\begin{align}
&\quad \frac{d}{dt}\Im \int_{\mathbb{R}^N}\bar{u}(x\cdot \nabla u)dx\nonumber\\
&=-2\int_{\mathbb{R}^N}|\nabla u|^2dx-(N+2)\int_{\mathbb{R}^N}|\nabla h(|u|^2)|^2dx\nonumber\\
&\quad-8N\int_{\mathbb{R}^N}h'(|u|^2)h''(|u|^2)|u|^4|\nabla u|^2dx+N\int_{\mathbb{R}^N}[|u|^2F(|u|^2)-G(|u|^2)]dx.\label{10131}
\end{align}
Lemma 2.1 is proved.\hfill $\Box$

{\bf Remark 2.1.} Although there doesn't exist the embrace relationship between the spaces $L^{p_1}(\mathbb{R}^N)$ and $L^{p_2}(\mathbb{R}^N)$ for $p_1>p_2>2$,
we can obtain the relationship between $\|u\|_{L^{p_1}}$ and $\|u\|_{L^{p_2}}$ if $u$ is the solution of (\ref{1}). In fact, using H\"{o}lder's inequality and the conservation law of mass, we have
\begin{align*}
\int_{\mathbb{R}^N}|u|^{p_2}dx&\leq \left(\int_{\mathbb{R}^N}|u|^2dx\right)^{\frac{p_1-p_2}{p_1-2}}\left(\int_{\mathbb{R}^N}|u|^{p_1}dx\right)^{\frac{p_2-2}{p_1-2}}\\
&=\left(\int_{\mathbb{R}^N}|u_0|^2dx\right)^{\frac{p_1-p_2}{p_1-2}}\left(\int_{\mathbb{R}^N}|u|^{p_1}dx\right)^{\frac{p_2-2}{p_1-2}}
\end{align*}

\section{The proofs of Theorem 1}

\qquad In this section, we will prove Theorem 1 and establish the sufficient conditions on the global existence of the solution to (\ref{1}).

{\bf Proof of Theorem 1:} Case (1). $F(s)=-F_2(s)\leq 0$ for $s\geq 0$. In this case $G(s)=-G_2(s)\leq 0$ for $s\geq 0$. The global existence of the solution is a direct result of the energy conversation law of Lemma 2.1(ii) because
$$
\int_{\mathbb{R}^N}|\nabla u|^2dx+\int_{\mathbb{R}^N}|\nabla h(|u|^2)|^2dx\nonumber\\
+\int_{\mathbb{R}^N}|G(|u|^2)|dx=2E(u_0),
$$
which implies that $\int_{\mathbb{R}^N}|\nabla u|^2dx+\int_{\mathbb{R}^N}|\nabla h(|u|^2)|^2dx+\int_{\mathbb{R}^N}|G(|u|^2)|dx$ is uniformly bounded for all $t>0$.

Case (2). $F(s)=F_1(s)\geq 0$ or $F(s)=F_1(s)$ changes sign for $s\geq 0$. $G(s)=G_1(s)$.
Denote
$$
\frac{1}{\tau_j}=\frac{q_j-1}{q_j-\theta_j},\qquad \frac{1}{\tau'_j}=\frac{1-\theta_j}{q_j-\theta_j},\quad j=1,2,
$$
Using the energy conversation law of Lemma 2.1(ii), using H\"{o}der's inequality,  Young's inequality, then Sobolev's inequality, after some elementary computation, we have
\begin{align}
&\quad\int_{\mathbb{R}^N}|\nabla u|^2dx+\int_{\mathbb{R}^N}|\nabla h(|u|^2)|^2dx\nonumber\\
&\leq C+\sum_{j=1}^2\left(c_j^{\frac{1}{\tau_j}}{c'}_j^{\frac{1}{\tau'_j}}\|u_0\|_2^2
+2^{\frac{2^*-1}{\tau'_j}}c_j^{\frac{1}{\tau_j}}(\epsilon_jC_s)^{\frac{1}{\tau'_j}}\|u_0\|_2^{\frac{2}{\tau_j}}
\left(\int_{\mathbb{R}^N}|\nabla u|^2dx\right)^{\frac{2^*}{2\tau'_j}}\right)\nonumber\\
&\quad+\sum_{j=1}^22^{\frac{2^*-1}{\tau'_j}}c_j^{\frac{1}{\tau_j}}(\epsilon_jC_s)^{\frac{1}{\tau'_j}}\|u_0\|_2^{\frac{2}{\tau_j}}\left(\int_{\mathbb{R}^N} |\nabla h(|u|^2)|^2dx\right)^{\frac{2^*}{2\tau'_j}}.\label{10151}
\end{align}

Now we discuss (\ref{10151}) in four subcases.

Subcase (a) $(2^*-2)\theta_1+2q_1>2^*$, $(2^*-2)\theta_2+2q_2>2^*$. Using (\ref{10151}), we have
\begin{align}
&\quad\int_{\mathbb{R}^N}|\nabla u|^2dx+\int_{\mathbb{R}^N}|\nabla h(|u|^2)|^2dx\nonumber\\
& \leq C(c_1,c_2,c'_1, c'_2,\epsilon_1,\epsilon_2, C_s,q_1, q_2, \theta_1,\theta_2,u_0)+\frac{1}{2}\int_{\mathbb{R}^N}[|\nabla u|^2dx+ |\nabla h(|u|^2)|^2]dx,
\end{align}
which implies that
$$\int_{\mathbb{R}^N}|\nabla u|^2dx+\int_{\mathbb{R}^N}|\nabla h(|u|^2)|^2dx\leq 2C(c_1,c_2,c'_1, c'_2, \epsilon_1,\epsilon_2, C_s,q_1, q_2, \theta_1,\theta_2,u_0),$$
where $C(c_1,c_2,c'_1, c'_2, \epsilon_1,\epsilon_2, C_s, q_1, q_2, \theta_1,\theta_2,u_0)$ is a positive constant depends on $c_1,c_2$, $c'_1, c'_2$, $\epsilon_1,\epsilon_2$, $C_s, q_1, q_2$, $\theta_1,\theta_2$ and $u_0$.

Subcase (b) $(2^*-2)\theta_1+2q_1=2^*$, $(2^*-2)\theta_2+2q_2>2^*$. Using (\ref{10151}), we obtain
\begin{align}
&\quad\int_{\mathbb{R}^N}|\nabla u|^2dx+\int_{\mathbb{R}^N}|\nabla h(|u|^2)|^2dx\nonumber\\
&\leq C'+\sum_{j=1}^2c_j^{\frac{1}{\tau_j}}{c'}_j^{\frac{1}{\tau'_j}}\|u_0\|_2^2\nonumber\\
&\quad+2^{\frac{2^*-1}{\tau'_1}}c_1^{\frac{1}{\tau_1}}(\epsilon_1C_s)^{\frac{1}{\tau'_1}}\|u_0\|_2^{\frac{2}{\tau_1}}
\left(\int_{\mathbb{R}^N}|\nabla u|^2dx+\int_{\mathbb{R}^N}|\nabla h(|u|^2)|^2dx\right)\nonumber\\
&\quad+\frac{1}{4}[1-2^{\frac{2^*-1}{\tau'_1}}c_1^{\frac{1}{\tau_1}}(\epsilon_1C_s)^{\frac{1}{\tau'_1}}\|u_0\|_2^{\frac{2}{\tau_1}}]
\left(\int_{\mathbb{R}^N}|\nabla u|^2dx+\int_{\mathbb{R}^N}|\nabla h(|u|^2)|^2dx\right).\label{6161}
\end{align}
Consequently,
\begin{align}
&\quad\frac{3}{4}[1-2^{\frac{2^*-1}{\tau'_1}}c_1^{\frac{1}{\tau_1}}(\epsilon_1C_s)^{\frac{1}{\tau'_1}}\|u_0\|_2^{\frac{2}{\tau_1}}]
\left(\int_{\mathbb{R}^N}|\nabla u|^2dx+\int_{\mathbb{R}^N}|\nabla h(|u|^2)|^2dx\right)\nonumber\\
&\leq \sum_{j=1}^2c_j^{\frac{1}{\tau_j}}{c'}_j^{\frac{1}{\tau'_j}}\|u_0\|_2^2+C'.\label{6162}
\end{align}

Subcase (c). $(2^*-2)\theta_1+2q_1>2^*$, $(2^*-2)\theta_2+2q_2=2^*$. Similar to Subcase (b), we can get
\begin{align}
&\quad\frac{3}{4}[1-2^{\frac{2^*-1}{\tau'_2}}c_2^{\frac{1}{\tau_2}}(\epsilon_2C_s)^{\frac{1}{\tau'_2}}\|u_0\|_2^{\frac{2}{\tau_2}}]
\left(\int_{\mathbb{R}^N}|\nabla u|^2dx+\int_{\mathbb{R}^N}|\nabla h(|u|^2)|^2dx\right)\nonumber\\
&\leq \sum_{j=1}^2c_j^{\frac{1}{\tau_j}}{c'}_j^{\frac{1}{\tau'_j}}\|u_0\|_2^2+C'.\label{6164}
\end{align}

Subcase (d) $(2^*-2)\theta_1+2q_1=2^*$, $(2^*-2)\theta_2+2q_2=2^*$. Similar to Subcase (b) and Subcase (c), we can get
\begin{align}
&\quad[1-\sum_{j=1}^22^{\frac{2^*-1}{\tau'_j}}c_j^{\frac{1}{\tau_j}}(\epsilon_jC_s)^{\frac{1}{\tau'_j}}\|u_0\|_2^{\frac{2}{\tau_j}}]
\left(\int_{\mathbb{R}^N}|\nabla u|^2dx+\int_{\mathbb{R}^N}|\nabla h(|u|^2)|^2dx\right)\nonumber\\
&\leq \sum_{j=1}^2c_j^{\frac{1}{\tau_j}}{c'}_j^{\frac{1}{\tau'_j}}\|u_0\|_2^2+C. \label{6166}
\end{align}
So we can know that $\int_{\mathbb{R}^N}|u|^2dx+\int_{\mathbb{R}^N}|\nabla u|^2dx+\int_{\mathbb{R}^N}|\nabla h(|u|^2)|^2dx$ is uniformly bounded for all $t>0$ in Case (2).

Case (3). $F(s)=F_1(s)-F_2(s)$.

Subcase (iv) $F_1(s)$ satisfies the assumptions of that in Case 2. Recalling
\begin{align*}
&\quad\int_{\mathbb{R}^N}|\nabla u|^2dx+\int_{\mathbb{R}^N}|\nabla h(|u|^2)|^2dx+\int_{\mathbb{R}^N}G_2(|u|^2)dx\nonumber\\
&=2E(u_0)+\int_{\mathbb{R}^N}G_1(|u|^2)dx,
\end{align*}
repeat the courses in Case (2), we can prove that $\int_{\mathbb{R}^N}|u|^2dx+\int_{\mathbb{R}^N}|\nabla u|^2dx+\int_{\mathbb{R}^N}|\nabla h(|u|^2)|^2dx$ is uniformly bounded for $t>0$.

Subcase (v) $|G_1(s)|\leq cs+ G_2(s)$ for $s\geq 0$.
Recalling
\begin{align*}
&\quad\int_{\mathbb{R}^N}|\nabla u|^2dx+\int_{\mathbb{R}^N}|\nabla h(|u|^2)|^2dx+\int_{\mathbb{R}^N}G_2(|u|^2)dx\nonumber\\
&\leq 2E(u_0)+c\int_{\mathbb{R}^N}|u|^2dx+\int_{\mathbb{R}^N}G_2(|u|^2)dx,
\end{align*}
which implies that $\int_{\mathbb{R}^N}|u|^2dx+\int_{\mathbb{R}^N}|\nabla u|^2dx+\int_{\mathbb{R}^N}|\nabla h(|u|^2)|^2dx$ is uniformly bounded for $t>0$.

Subcase (vi) $F_1(s)$ and $F_2(s)$ satisfy (\ref{6171}) and (\ref{6172}). Using Young inequality,  we have
\begin{align}
&\quad\int_{\mathbb{R}^N}|\nabla u|^2dx+\int_{\mathbb{R}^N}|\nabla h(|u|^2)|^2dx+\int_{\mathbb{R}^N}G_2(|u|^2)dx\nonumber\\
&=2E(u_0)+\int_{\mathbb{R}^N}G_1(|u|^2)dx\nonumber\\
&\leq 2E(u_0)+\int_{\{|u|\leq 1\}}|G_1(|u|^2)|dx+\int_{\{|u|>1\}}|G_1(|u|^2)|dx\nonumber\\
&\leq 2E(u_0)+\int_{\{|u|\leq 1\}}\left(C(\tilde{\alpha}_1,\tilde{\beta}_1,\tilde{c}_1,\tilde{c}_2)[|G_1(|u|^2)|]^{\tilde{\alpha}_1}
+\frac{1}{4\tilde{c}_2}[|G_1(|u|^2)|]^{\tilde{\beta_1}}\right)dx\nonumber\\
&\quad+\int_{\{|u|>1\}}\left(C(\tilde{\alpha}_2,\tilde{\beta}_2,\tilde{c}'_1,\tilde{c}'_2)[|G_1(|u|^2)|]^{\tilde{\alpha}_2}
+\frac{1}{4\tilde{c}'_2}[|G_1(|u|^2)|]^{\tilde{\beta}_2}\right)dx\nonumber\\
&\leq 2E(u_0)+C\int_{\mathbb{R}^N}|u_0|^2dx+\frac{1}{2}\int_{\mathbb{R}^N}G_2(|u|^2)dx,
\end{align}
which implies that $\int_{\mathbb{R}^N}|u|^2dx+\int_{\mathbb{R}^N}|\nabla u|^2dx+\int_{\mathbb{R}^N}|\nabla h(|u|^2)|^2dx$ is uniformly bounded for all $t>0$.\hfill $\Box$

Noticing that $\|u(\cdot,t)\|_{L^2}=\|u_0\|_{L^2}$, using the results of Theorem 1, we can get some related results below.

{\bf Proposition 3.1.} {\it Assume that, excluding (v) in Case (3), the other conditions of Theorem 1 hold and $u$ is the global solution of (\ref{1}).
Suppose that the functions $f(s)$ and $g(s)$ satisfying the following conditions: There exist $0<\alpha_1<1$, $0<\alpha_2<1$, $\beta_1>1$ and $\beta_2>1$ such that
$$
(EC)\quad [|f(s)|]^{\alpha_1}\leq c_1s ,\ [|f(s)|]^{\beta_1}\leq C_1[h(s)]^{2^*},\ [|g(s)|]^{\alpha_2}\leq c_2s,\ [|g(s)|]^{\beta_2}\leq C_2|G(s)|
$$
for $s\geq 0$, where $c_1, c_2, C_1$ and $C_2$ are positive constants. Let
\begin{align}
\tau_1=\frac{(\beta_1-\alpha_1)}{(1-\alpha_1)},\quad \tau_2=\frac{(\beta_2-\alpha_2)}{(1-\alpha_2)}.\label{6110}
\end{align}
Then
\begin{align}
\int_{\mathbb{R}^N}|f(|u|^2)|dx\leq C,\quad \int_{\mathbb{R}^N}|g(|u|^2)|dx\leq C.\label{6111}
\end{align}
}

{\bf Proof:} Noticing that
\begin{align}
\int_{\mathbb{R}^N}|f(|u|^2)|dx&\leq \left(\int_{\mathbb{R}^N}[|f(|u|^2)|]^{\alpha_1}dx\right)^{\frac{1}{\tau'_1}}
\left(\int_{\mathbb{R}^N}[|f(|u|^2)|]^{\beta_1}dx\right)^{\frac{1}{\tau_1}}\nonumber\\
&\leq C\left(\int_{\mathbb{R}^N}|u|^2dx\right)^{\frac{1}{\tau'_1}}
\left(\int_{\mathbb{R}^N}[h(|u|^2)]^{2^*}dx\right)^{\frac{1}{\tau_1}}\nonumber\\
&\leq C\|u_0\|_2^{\frac{2}{\tau'_1}}\left(\int_{\mathbb{R}^N}|\nabla h(|u|^2)|^2dx\right)^{\frac{2^*}{2\tau_1}}\label{6115}
\end{align}
and
\begin{align}
\int_{\mathbb{R}^N}|g(|u|^2)|dx&\leq \left(\int_{\mathbb{R}^N}[|g(|u|^2)|]^{\alpha_2}dx\right)^{\frac{1}{\tau'_2}}
\left(\int_{\mathbb{R}^N}[|g(|u|^2)|]^{\beta_2}dx\right)^{\frac{1}{\tau_2}}\nonumber\\
&\leq C\left(\int_{\mathbb{R}^N}|u|^2dx\right)^{\frac{1}{\tau'_2}}
\left(\int_{\mathbb{R}^N}|G(|u|^2)|dx\right)^{\frac{1}{\tau_2}}\nonumber\\
&\leq C\|u_0\|_2^{\frac{2}{\tau'_1}}\left(\int_{\mathbb{R}^N}|G(|u|^2)|dx\right)^{\frac{1}{\tau_2}},\label{6116}
\end{align}
we can obtain the conclusions.\hfill $\Box$

We would like to give an example to illustrate the conclusions of Proposition 3.1.

$h(s)=as^{p}(a>0, p>\frac{1}{2^*})$ and $F(s)=bs^{\tilde{q}}(b<0, \tilde{q}>0)$, that is, $h(|u|^2)=a|u|^{2p}(a>0, p>\frac{1}{2^*})$ and
$|G(|u|^2)|=\frac{-b}{\tilde{q}+1}|u|^{2\tilde{q}+2}(b<0,\tilde{q}>0)$. Then for any $1<r_1\leq p\cdot 2^*$ and $1<r_2\leq\tilde{q}+1$,
we have

\begin{align}
\|u\|_{L^{2r_1}}^{2r_1}&=\int_{\mathbb{R}^N}|u|^{2r_1}dx\leq \left(\int_{\mathbb{R}^N}|u|^2dx\right)^{\frac{1}{\tau'_1}}
\left(\int_{\mathbb{R}^N}|u|^{2p\cdot 2^*}dx\right)^{\frac{1}{\tau_1}}\nonumber\\
&\leq C\left(\int_{\mathbb{R}^N}|u_0|^2dx\right)^{\frac{1}{\tau'_1}}
\left(\int_{\mathbb{R}^N}|\nabla (|u|^{2p})|^2dx\right)^{\frac{2^*}{2\tau_1}}\nonumber\\
&\leq C,\label{6117}\\
\|u\|_{L^{2r_2}}^{2r_2}&=\int_{\mathbb{R}^N}|u|^{2r_2}dx\leq \left(\int_{\mathbb{R}^N}|u|^2dx\right)^{\frac{1}{\tau'_2}}
\left(\int_{\mathbb{R}^N}|u|^{2\tilde{q}+2}dx\right)^{\frac{1}{\tau_2}}\nonumber\\
&\leq C\left(\int_{\mathbb{R}^N}|u_0|^2dx\right)^{\frac{1}{\tau'_1}}
\left(\int_{\mathbb{R}^N}|u|^{2\tilde{q}+2}dx\right)^{\frac{1}{\tau_2}}\nonumber\\
&\leq C.\label{6118}
\end{align}

Similarly, recalling that $\|\nabla u(\cdot,t)\|_{L^2}\leq C$ uniformly for all $t>0$, we have the following propositions.

{\bf Proposition 3.2.} {\it  Assume that, excluding (v) in Case (3), the conditions of Theorem 1 hold and $u$ is the global solution of (\ref{1}).
Suppose that $f(s)$ satisfying
$$\left([f'(s)]^2s\right)^{\tilde{\tau}}\leq C[h'(s)]^2s$$
for $s\geq 0$, where the constant $\tilde{\tau}>1$, and $g(s)$ satisfying
$$[g'(s)]^2s\leq C |G_1(s)|\quad  {\rm or} \quad [g'(s)]^2s\leq C G_2(s).$$
Then
\begin{align}
\|\nabla f(|u|^2)\|_{L^2}^2=\int_{\mathbb{R}^N}|\nabla f(|u|^2)|^2dx\leq C, \quad \|\nabla g(|u|^2)\|_{L^1}=\int_{\mathbb{R}^N}|\nabla g(|u|^2)|dx\leq C.\label{6121}
\end{align}
}

{\bf Proof:} Noticing that
\begin{align*}
\int_{\mathbb{R}^N}|\nabla f(|u|^2)|^2dx&\leq \left(\int_{\mathbb{R}^N}|\nabla u|^2dx\right)^{\frac{1}{\tilde{\tau}'}}
\left(\int_{\mathbb{R}^N}\left([f'(|u|^2)]^2|u|^2\right)^{\tilde{\tau}}|\nabla u|^2dx\right)^{\frac{1}{\tilde{\tau}}}\\
&\leq C\left(\int_{\mathbb{R}^N}|\nabla h(|u|^2)|^2dx\right)^{\frac{1}{\tilde{\tau}}}\leq C'.
\end{align*}

We only prove the case of $[g'(s)]^2s\leq G_1(s)$, the proof of the other case is similar. Noticing that
\begin{align*}
\int_{\mathbb{R}^N}|\nabla g(|u|^2)|dx&\leq \left(\int_{\mathbb{R}^N}|\nabla u|^2dx\right)^{\frac{1}{2}}
\left(\int_{\mathbb{R}^N}[4g'(|u|^2)]^2|u|^2dx\right)^{\frac{1}{2}}\\
&\leq C\left(\int_{\mathbb{R}^N}|G_1(|u|^2)|dx\right)^{\frac{1}{2}}\leq C',
\end{align*}
we can get the conclusions.\hfill $\Box$

\section{The Proof of Theorem 2}
\qquad In this section, we will give the proof of Theorem 2 and deal with the sufficient conditions on blowup in finite time for the solution by using the results of Lemma 2.1.

{\bf Proof of Theorem 2:} Wherever $u$ exists, let
$$
y(t)=\Im \int_{\mathbb{R}^N}\bar{u}(x\cdot \nabla u)dx.
$$

We discuss it in three cases: Case 1. $h(s)\equiv 0$;  Case 2. $h(s)\neq 0$ and $k\leq -\frac{1}{2}$;  Case 3. $h(s)\neq 0$ and $k>-\frac{1}{2}$.

First, we deal with it in Case 2. By the result (iv) of Lemma 2.1 and the assumption of $$2(c_N-1)E(u_0)+c_MM(u_0)\leq 0\qquad {\rm if} \quad k\leq-\frac{1}{2},$$ we have
\begin{align}
\dot{y}(t)&\geq [N(c_N-1)-2]\int_{\mathbb{R}^N}|\nabla u|^2dx.\label{10141}
\end{align}

Case 1. Similar to the computations above, (\ref{10141}) is still hold in this case.

Case 3. By the result (iv) of Lemma 2.1 and the assumption of $$2[(2k+1)N+2]E(u_0)+c_MM(u_0)\leq 0\qquad {\rm if} \quad k>-\frac{1}{2},$$ we have

\begin{align}
\dot{y}(t)&\geq (2k+1)N\int_{\mathbb{R}^N}|\nabla u|^2dx.\label{10141'}
\end{align}

In a word, $\dot{y}(t)\geq 0$ in the three cases above.

Under the conditions of $y(0)=\Im \int_{\mathbb{R}^N}\bar{u}_0(x\cdot u_0)dx>0$, we know that $y(t)$ is increasing, which implies that $y(t)>0$ wherever $u$ exists.

On the other hand, by the result (iii) of Lemma 2.1,
$$
\frac{d}{dt}\int_{\mathbb{R}^N}|x|^2|u|^2dx
=-4\Im \int_{\mathbb{R}^N}\bar{u}(x\cdot \nabla u)dx=-4y(t)<0,
$$
which means that
$$
\int_{\mathbb{R}^N}|x|^2|u|^2dx\leq \int_{\mathbb{R}^N}|x|^2|u_0|^2dx:=d^2_0<+\infty.
$$
Using Schwarz inequality, we get
\begin{align}
y(t)\leq \left(\int_{\mathbb{R}^N}|x|^2|u|^2dx\right)^{\frac{1}{2}}\left(\int_{\mathbb{R}^N}|\nabla u|^2dx\right)^{\frac{1}{2}}
\leq d_0\left(\int_{\mathbb{R}^N}|\nabla u|^2dx\right)^{\frac{1}{2}}.\label{10142}
\end{align}
(\ref{10141}) and (\ref{10142}) imply that
\begin{align*}
\dot{y}(t)\geq \frac{[N(c_N-1)-2]}{d_0^2}y^2(t),
\end{align*}
while (\ref{10141'}) and (\ref{10142}) imply that
\begin{align*}
\dot{y}(t)\geq \frac{(2k+1)N}{d_0^2}y^2(t)
\end{align*}
Then we have
\begin{align}
\dot{y}(t)\geq \frac{\max([N(c_N-1)-2],(2k+1)N)}{d_0^2}y^2(t):=\frac{c(k,N)}{d_0^2}y^2(t)
\label{10143}
\end{align}
with $y(0)>0$. Integrating (\ref{10143}), we obtain
$$
y(t)\geq \frac{y(0)d^2_0}{d_0^2-y(0)c(k,N)t}, \quad 0\leq t<\frac{d_0^2}{y(0)c(k,N)}.
$$
Consequently,
$$
\|\nabla u\|_2\geq \frac{y(0)d_0}{d_0^2-y(0)c(k,N)t},
$$
and there exist $T\leq T_0=\frac{d_0^2}{y(0)c(k,N)}$ such that
$$
\lim_{t\rightarrow T^-} \|\nabla u\|_2=+\infty.
$$
Theorem 2 is proved. \hfill $\Box$

As a direct result of Theorem 1 and Theorem 2, we will give a corollary below and compare our results with those of others.

{\bf Corollary 4.1.} {\it Assume that $u(x,t)$ is the solution to
\begin{equation}
\label{3'} \left\{
\begin{array}{lll}
iu_t=\Delta u+2\alpha(\Delta |u|^{2\alpha})|u|^{2\alpha-2}u+|u|^{q-2}u=0\quad {\rm for} \ x\in \mathbb{R}^N, \ t>0\\
u(x,0)=u_0(x),\quad x\in \mathbb{R}^N.
\end{array}\right.
\end{equation}
Here $\alpha>0$, $q>2$, $u_0\in H^1(\mathbb{R}^N)$. Then $u(x,t)$ is global existence if $2<q<4\max\{\alpha,\frac{1}{2}\}+\frac{4}{N}$ for any $u_0\in H^1(\mathbb{R}^N)$, while it will blow up in finite time
if $4\max\{\alpha,\frac{1}{2}\}+\frac{4}{N}\leq q<2\max\{\alpha,\frac{1}{2}\}\cdot2^*$ for $u_0$ satisfying (i) $E(u_0)\leq 0$;  (ii) $\Im \int_{\mathbb{R}^N}\bar{u}_0(x\cdot u_0)dx>0$, $|x|u_0\in L^2(\mathbb{R}^N)$. }

{\bf Proof:} Letting $h(s)=s^{\alpha}$, $F(s)=s^{\frac{q-2}{2}}$, taking $\theta_1=\theta_2=\frac{2}{q}$, $q_1=\frac{q2^*+2}{q}$,
$q_2=\frac{2^*}{q}\max(2\alpha,1)$  in Theorem 1 and $c_M=0$ and $c_N=\frac{q}{2}$ in Theorem 2, we can get the conclusions.\hfill $\Box$

\section{The proof of Theorem 3}

\qquad In this section, we will prove Theorem 3 and establish a sharp threshold for the blowup and global existence of the solution to (\ref{1}) under certain conditions.

{\bf The proof of Theorem 3. } We proceed in four Steps.

Step 1. We will prove $d_I>0$.

Since $Q(w)=0$, $w\not\equiv 0$, after some elementary computations,  we have
\begin{align}
&\quad \underline{l}(\int_{\mathbb{R}^N}|\nabla w|^2dx+\int_{\mathbb{R}^N}|\nabla h(|w|^2)|^2dx)\nonumber\\
&\leq N \sum_{j=3}^4\left(\int_{\mathbb{R}^N}c_j|w|^2dx\right)^{\frac{1}{\tau_j}}2^{2*}(c'_jC_s)^{\frac{1}{\tau'_j}}
\left(\int_{\mathbb{R}^N}([|\nabla w|^2+|\nabla h(|w|^2)|^2]dx\right)^{\frac{2^*}{2\tau'_j}},\label{6291}
\end{align}
where
$$
 \frac{1}{\tau_j}=\frac{q_j-1}{q_j-\theta_j},\qquad \frac{1}{\tau'_j}=\frac{1-\theta_j}{q_j-\theta_j},\quad j=3,4.
$$

If $\frac{2^*}{2\tau'_j}=1$, i.e., $(2^*-2)\theta_j+2q_j=2^*(j=3,4)$, we get $\int_{\mathbb{R}^N}|w|^2dx\geq C$.

If $\frac{2^*}{2\tau'_j}>1$, i.e., $(2^*-2)\theta_j+2q_j=2^*(j=3,4)$, then (\ref{6291}) implies that
\begin{align}
\underline{l}\leq N\sum_{j=3}^42^{2^*}c_j^{\frac{1}{\tau_j}}(c'_jC_s)^{\frac{1}{\tau'_j}}\left(\int_{\mathbb{R}^N}|w|^2dx\right)^{\frac{1}{\tau_j}}
\left(\int_{\mathbb{R}^N}([|\nabla w|^2+|\nabla h(|w|^2)|^2]dx\right)^{\frac{2^*}{2\tau'_j}-1}.\label{1018w1}
\end{align}
Using Young's inequality to (\ref{1018w1}), we have
\begin{align}
\underline{l}\leq N \sum_{j=3}^42^{2^*}c_j^{\frac{1}{\tau_j}}(c'_jC_s)^{\frac{1}{\tau'_j}}\left\{\int_{\mathbb{R}^N}|w|^2dx+
\left(\int_{\mathbb{R}^N}[|\nabla w|^2+|\nabla h(|w|^2)|^2]dx\right)^{(\frac{2^*}{2\tau'_j}-1)\tau'_j}\right\},\label{1018w2}
\end{align}
which implies that
\begin{align}
\int_{\mathbb{R}^N}|w|^2dx+\int_{\mathbb{R}^N}|\nabla w|^2dx+\int_{\mathbb{R}^N}|\nabla h(|w|^2)|^2dx\geq C>0\label{1018w3}
\end{align}
if $Q(w)=0$.

If one of $\frac{2^*}{2\tau'_3}$ and $\frac{2^*}{2\tau'_4}$ equals 1, while another is larger than 1, we can using Young's inequality and similarly
deal with (\ref{1018w1}) to get (\ref{1018w3}).

On the other hand, using $Q(w)=0$ again, we get
\begin{align}
&\quad L(\int_{\mathbb{R}^N}|\nabla w|^2dx+\int_{\mathbb{R}^N}|\nabla h(|w|^2)|^2dx)\nonumber\\
&=N\int_{\mathbb{R}^N}[|w|^2F(|w|^2)-G(|w|^2)]dx\geq N[c_N-1)\int_{\mathbb{R}^N}G(|w|^2)dx.\label{1019s1}
\end{align}
Therefore
\begin{align}
E(w)&=\frac{1}{2}\left(\int_{\mathbb{R}^N}|\nabla w|^2dx+\int_{\mathbb{R}^N}|\nabla h(|w|^2)|^2dx-\int_{\mathbb{R}^N}G(|w|^2)dx\right)\nonumber\\
&\geq \frac{1}{2}\left(1-\frac{L}{N(c_N-1)}\right)\left(\int_{\mathbb{R}^N}|\nabla w|^2dx+\int_{\mathbb{R}^N}|\nabla h(|w|^2)|^2dx\right).\label{1018w4}
\end{align}
(\ref{1019s1}) and (\ref{1018w4}) mean that
\begin{align}
&\quad\frac{\omega}{2} \int_{\mathbb{R}^N}|w|^2dx+E(w)\nonumber\\
&\geq \frac{1}{2}\min\left(\omega, 1-\frac{L}{N(c_N-1)}\right)\left( \int_{\mathbb{R}^N}|w|^2dx+\int_{\mathbb{R}^N}|\nabla w^2dx+\int_{\mathbb{R}^N}|\nabla h(|w|^2)|^2dx\right)\nonumber\\
&\geq C>0.\label{1018w3'}
\end{align}
Therefore $d_I>0$.

Step 2. Denote
$$
K_+=\{u\in H^1(\mathbb{R}^N)\setminus\{0\},\ Q(u)>0,\
\frac{\omega}{2}\|u\|_2^2+E(u)<d_I\}
$$
and
$$
K_-=\{u\in H^1(\mathbb{R}^N)\setminus\{0\},\ Q(u)<0,\
\frac{\omega}{2}\|u\|_2^2+E(u)<d_I\}.
$$
We will prove that $K_+$ and $K_-$ are invariant sets of
(\ref{1}).

Assume that $u_0\in K_+$, i.e., $Q(u_0)>0$ and
$\frac{\omega}{2}\|u_0\|_2^2+E(u_0)<d_I$. It is easy to verify that
\begin{align}
u(\cdot,t)\in H^1(\mathbb{R}^N)\setminus\{0\}, \quad
\frac{\omega}{2}\|u(\cdot,t)\|_2^2+E(u(\cdot,t))<d_I\label{9626z2}\end{align}
because $\|u\|^2_2$ and $E(u)$ are conservation quantities for (\ref{1}).

We need to show that $Q(u(\cdot,t))>0$ for $t\in (0,T)$. Contradictorily, if there exists $t_1\in (0,T)$ such that
$Q(u(\cdot,t_1))<0$, then there exists a $t_2\in [0, t_1]$ such that
$Q(u(\cdot,t_2))=0$ by the continuity. And
$$\frac{\omega}{2}\|u(\cdot,t_2)\|_2^2+E(u(\cdot,t_2))<d_I$$
by (\ref{9626z2}), which is a contradiction to the definition of $d_I$.
Hence $Q(u(\cdot,t))>0$. This inequality and (\ref{9626z2}) imply
that $u(\cdot,t)\in K_+$, which means that $K_+$ is a invariant
set of (\ref{1}).

Similarly, we can prove that $K_-$ is also a
invariant set of (\ref{1}). We omit the details here.

Step 3. Assume that $Q(u_0)>0$ and
$\frac{\omega}{2}\|u_0\|_2^2+E(u_0)<d_I$. Since $\mathcal{K}$ is invariant set of (\ref{1}), we
have $Q(u(\cdot,t))>0$ and
$\frac{\omega}{2}\|u(\cdot,t)\|_2^2+E(u(\cdot,t))<d_I$.  Using $Q(u(\cdot,t))>0$, we get
\begin{align}
&\quad L(\int_{\mathbb{R}^N}|\nabla u(\cdot,t)|^2dx+\int_{\mathbb{R}^N}|\nabla h(|u(\cdot,t)|^2)|^2dx)\nonumber\\
&\geq N\int_{\mathbb{R}^N}[|u|^2F(|u(\cdot,t)|^2)-G(|u(\cdot,t)|^2)]dx\geq N(c_N-1)\int_{\mathbb{R}^N}G(|u(\cdot,t)|^2)dx.\label{1019s2}
\end{align}
Using (\ref{1019s2}), we obtain
\begin{align}
E(u(\cdot,t))&=\frac{1}{2}\int_{\mathbb{R}^N}|\nabla u(\cdot,t)|^2dx+\int_{\mathbb{R}^N}|\nabla h(|u(\cdot,t)|^2)|^2dx-\int_{\mathbb{R}^N}G(|u(\cdot,t)|^2)dx\nonumber\\
&\geq \frac{1}{2}\left(1-\frac{L}{N(c_N-1)}\right)\left(\int_{\mathbb{R}^N}|\nabla u(\cdot,t)|^2dx+\int_{\mathbb{R}^N}|\nabla h(|u(\cdot,t)|^2)|^2dx\right).\label{1019s3}
\end{align}
By mass conversation law $\|u(\cdot,t)\|_2^2=\|u_0\|_2^2$, using (\ref{1019s2}) and (\ref{1019s3}), we get
\begin{align*}
d_I&>\frac{\omega}{2} \|u(\cdot,t)\|_2^2+E(u(\cdot,t))\nonumber\\
&\geq \frac{1}{2}\left(1-\frac{L}{N(c_N-1)}\right)\left(\int_{\mathbb{R}^N}|\nabla u(\cdot,t)|^2dx+\int_{\mathbb{R}^N}|\nabla h(|u(\cdot,t)|^2)|^2dx\right)
\end{align*}
and
$$
\int_{\mathbb{R}^N}|\nabla u(\cdot,t)|^2dx+\int_{\mathbb{R}^N}|\nabla h(|u(\cdot,t)|^2)|^2dx\leq C<\infty,
$$
i.e., the solution $u(x,t)$ of (\ref{1}) exists globally.

Step 4. Suppose that $|x|u_0\in L^2(\mathbb{R}^N)$, $Q(u_0)<0$ and
$\frac{\omega}{2}\|u_0\|_2^2+E(u_0)<d_I$. Since $\mathcal{K}_-$ is a invariant set of (\ref{1}), we have $Q(u(\cdot,t))<0$ and
$\frac{\omega}{2}\|u(\cdot,t)\|_2^2+E(u(\cdot,t))<d_I$.

Let $J(t)=\int_{\mathbb{R}^N} |x|^2|u|^2dx$. Then
$$
J''(t)=4Q(u(x,t)),\quad J'(t)=-4\Im \int_{\mathbb{R}^N} \bar{u}(x\cdot \nabla u)dx.
$$
Since
$J'(0)=-4\Im \int_{\mathbb{R}^N} \bar{u}_0(x\cdot \nabla u_)dx<0$, we have
$$
J'(t)=J'(0)+\int_0^t J''(\tau)d\tau=J'(0)+4\int_0^t Q(u(\cdot,\tau))d\tau<J'(0)<0
$$
and
$$
J(t)=J(0)+\int_0^t J'(\tau)d\tau<J(0)+J'(0)t,
$$
which implies that the maximum existence interval for $t$ is finite, by the proof of Theorem 2,
we know that the solution blows up in finite time.\hfill
$\Box$

We also give a corollary of Theorem 3 as follows.

{\bf Corollary 5.1.} {\it Assume that $u(x,t)$ is the solution to (\ref{3'}), $\alpha\geq \frac{1}{2}$ and $4\alpha+\frac{4}{N}<q<2\alpha\cdot 2^*$. Moreover suppose that
there exists $\omega>0$ such that
\begin{align} d_I:=\inf_{\{w\in H^1(\mathbb{R}^N)\setminus \{0\};
Q(w)=0\}}\left(\frac{\omega}{2}\|w\|_2^2+E(w)\right)>0,\label{9651'}\end{align}
where
\begin{align}
Q(w)&=2\int_{\mathbb{R}^N}|\nabla w|^2dx+[(2\alpha-1)N+2]\int_{\mathbb{R}^N}|\nabla w^{2\alpha}|^2dx-\frac{(q-2)N}{q}\int_{\mathbb{R}^N}|w|^{q}dx,\\
E(w)&=\frac{1}{2}\int_{\mathbb{R}^N}[|\nabla w|^2+|\nabla w^{2\alpha}|^2]dx-\frac{1}{q}\int_{\mathbb{R}^N}|w|^{q}dx
\end{align}
and $u_0$
satisfies
$$\frac{\omega}{2}\|u_0\|_2^2+E(u_0)<d_I.$$

Then we have:

(1). If $Q(u_0)>0$, the solution of (\ref{1}) exists
globally;

(2). If $Q(u_0)<0$ and $|x|u_0\in L^2(\mathbb{R}^N)$, $\Im \int_{\mathbb{R}^N} \bar{u}_0(x\cdot \nabla u_0)dx>0$, the solution
of (\ref{1}) blows up in finite time.}

{\bf  Proof:} Letting $h(s)=s^{\alpha}$, $\alpha\geq \frac{1}{2}$, and $F(s)=s^{\frac{q-2}{2}}$, $q>2$, we can verify the conditions of Theorem 3, which implies the conclusions of this corollary are true.\hfill $\Box$

\section{The Pseudo-conformal Conservation Law and Asymptotic Behavior for the Global Solution}
\qquad In this section, we will  prove the pseudo-conformal conservation law and consider asymptotic behavior for the global solution of (\ref{1}).

{\bf Proof of Theorem 4:} 1. Assume that $u$ is the global solution of (\ref{1}),  $u_0\in X$ and $xu_0\in L^2(\mathbb{R}^N)$. Since
$$
E(u)=\frac{1}{2}\int_{\mathbb{R}^N}[|\nabla u|^2+|\nabla h(|u|^2)|^2-G(|u|^2)]dx
=E(u_0),
$$
we have
\begin{align}
P(t)&=\int_{\mathbb{R}^N}|xu|^2dx+4t\Im \int_{\mathbb{R}^N}\bar{u}(x\cdot \nabla u)dx+4t^2\int_{\mathbb{R}^N}|\nabla u|^2dx\nonumber\\
&\qquad+4t^2\int_{\mathbb{R}^N}|\nabla h(|u|^2)|^2dx-4t^2\int_{\mathbb{R}^N}G(|u|^2)dx\nonumber\\
&=\int_{\mathbb{R}^N}|xu|^2dx+4t\Im \int_{\mathbb{R}^N}\bar{u}(x\cdot \nabla u)dx+8t^2E(u_0).\label{692}
\end{align}
Recalling that
$$\frac{d}{dt} \int_{\mathbb{R}^N}|x|^2|u|^2dx=-4\Im \int_{\mathbb{R}^N} \bar{u}(x\cdot \nabla u)dx,$$
we get
\begin{align}
P'(t)&=4t\int_{\mathbb{R}^N}-4N[2h''(|u|^2)h'(|u|^2)|u|^2+(h'(|u|^2))^2]|u|^2|\nabla u|^2dx\nonumber\\
&\quad+ 4t\int_{\mathbb{R}^N}[N|u|^2F(|u|^2)-(N+2)G(|u|^2)]dx.\label{693}
\end{align}
Integrating (\ref{693}) from $0$ to $t$, we have
$$
P(t)=P(0)+4\int_0^t\tau\theta(\tau)d\tau=\int_{\mathbb{R}^N}|xu_0|^2dx+4\int_0^t\tau\theta(\tau)d\tau.
$$
That is,
\begin{align}
&\quad \int_{\mathbb{R}^N}|(x-2it\nabla)u|^2dx+4t^2\int_{\mathbb{R}^N}|\nabla h(|u|^2)|^2dx-4t^2\int_{\mathbb{R}^N}G(|u|^2)dx\nonumber\\
&=\int_{\mathbb{R}^N}|xu_0|^2dx+4\int_0^t\tau\theta(\tau)d\tau,\label{694}
\end{align}
where $\theta(\tau)$ is defined by (\ref{691'}).

2. Assume that $u$ is the blowup solution of (\ref{1}),  $u_0\in X$ and $xu_0\in L^2(\mathbb{R}^N)$. Using $E(u)=E(u_0)$, we have
\begin{align}
B(t)&:=\int_{\mathbb{R}^N}|(x+2i(T-t)\nabla)u|^2dx+4(T-t)^2\int_{\mathbb{R}^N}|\nabla h(|u|^2)|^2dx-4(T-t)^2\int_{\mathbb{R}^N}G(|u|^2)dx\nonumber\\
&=\int_{\mathbb{R}^N}|xu|^2dx-4(T-t)\Im \int_{\mathbb{R}^N}\bar{u}(x\cdot \nabla u)dx+8(T-t)^2E(u_0)\label{891}
\end{align}
and
\begin{align}
B'(t)
&=4(T-t)\left\{\int_{\mathbb{R}^N}4N[2h''(|u|^2)h'(|u|^2)|u|^2+( h'(|u|^2))^2]|u|^2|\nabla u|^2dx\right.\nonumber\\
&\qquad\left.+\int_{\mathbb{R}^N} [(N+2)G(|u|^2)-NF(|u|^2)|u|^2]dx\right\}.\label{892}
\end{align}
Integrating (\ref{892}) from $0$ to $t$, we have
\begin{align*}
B(t)&=B(0)-4\int_0^t(T-\tau)\theta(\tau)d\tau\nonumber\\
&=\int_{\mathbb{R}^N}|(x+2iT\nabla)u_0|^2dx+4T^2\int_{\mathbb{R}^N}[|\nabla h(|u_0|^2)|^2-G(|u_0|^2)]dx\nonumber\\
&\quad -4\int_0^t(T-\tau)\theta(\tau)d\tau,
\end{align*}
where $\theta(\tau)$ is defined by (\ref{691'}).\hfill $\Box$

Using Theorem 4 we derive the following result on asymptotic behaviors of the solution to (\ref{1}).

{\bf Theorem 5.} {\it 1. Assume that $u$ is the global solution of (\ref{1}),  $u_0\in X$ and $xu_0\in L^2(\mathbb{R}^N)$. Suppose that $F(s)=F_1(s)-F_2(s)$, $F_1(s)\geq 0$ or $F_1(s)$ changes sign for $s\geq 0$, $F_2(s)\geq 0$ for $s\geq 0$, and there exist $c_1, c'_1, c_2, c'_2>0$, $0<\gamma_1, \tilde{\gamma}_1<1$ and $\gamma_2,\tilde{\gamma}_2>1$ such that
\begin{align}
&\frac{2^*(1-\gamma_1)}{2(\gamma_2-\gamma_1)}=1,\quad \frac{2^*(1-\tilde{\gamma}_1)}{2(\tilde{\gamma}_2-\tilde{\gamma}_1)}=1,\label{6261}\\
& [|G_1(s)|]^{\gamma_1}\leq c_1s ,\quad [|G_1(s)|]^{\gamma_2}\leq c'_1[h(s)]^{2^*}\ {\rm for}\ 0\leq s\leq 1,\label{6262}\\
& [|G_1(s)|]^{\tilde{\gamma}_1}\leq c_2s ,\quad [|G_1(s)|]^{\tilde{\gamma}_2}\leq c'_2[h(s)]^{2^*}\ {\rm for}\  s>1.\label{6263}
\end{align}
Moreover, assume that
\begin{align}
\sum_{j=1}^2(c_j\|u_0\|^2_{L^2})^{\frac{1}{\tilde{\tau}'_j}}(c'_jC_s)^{\frac{1}{\tilde{\tau}_j}}<1.\label{61210}
\end{align}
Here
\begin{align}
\frac{1}{\tilde{\tau}'_1} = \frac{(\gamma_2-1)}{(\gamma_2-\gamma_1)},\quad \frac{1}{\tilde{\tau}_1} = \frac{(1-\gamma_1)}{(\gamma_2-\gamma_1)},\label{6129}\\
\frac{1}{\tilde{\tau}'_2} = \frac{(\tilde{\gamma}_2-1)}{(\tilde{\gamma}_2-\tilde{\gamma}_1)},\quad \frac{1}{\tilde{\tau}_2} = \frac{(1-\tilde{\gamma}_1)}{(\tilde{\gamma}_2-\tilde{\gamma}_1)}.\label{6130}
\end{align}

Then the following properties hold:

(1) If $2h''(s)h'(s)s+(h'(s))^2\geq 0$, $NF_1(s)s-(N+2)G_1(s)\leq 0$ and $NF_2(s)s-(N+2)G_2(s)\geq0 $, then there exists $C$ such that
\begin{align}
\int_{\mathbb{R}^N}[|\nabla h(|u|^2)|^2+|G_1(|u|^2)|+G_2(|u|^2)]dx\leq Ct^{-2}\quad {\rm for}\ t\geq 1. \label{6101}
\end{align}

(2) If $2h''(s)h'(s)s+(h'(s))^2\geq 0$, $NF_1(s)s-(N+2)G_1(s)\leq 0$ and $-k_1G_2(s)\leq NF_2(s)s-(N+2)G_2(s)<0$ for some $0<k_1<2$, then there exists $C$ such that
\begin{align}
\int_{\mathbb{R}^N}[|\nabla h(|u|^2)|^2+|G_1(|u|^2)|+G_2(|u|^2)]dx\leq\frac{C}{t^{2-k_1}}\quad {\rm for}\ t\geq 1.\label{6102}
\end{align}

(3) If $2h''(s)h'(s)s+(h'(s))^2\geq 0$, $0<NF_1(s)s-(N+2)G_1(s)\leq\tilde{k}_1|G_1(s)|$ for some \begin{align}
0<\tilde{k}_1<\frac{2[1-\sum_{j=1}^2(c_j\|u_0\|^2_{L^2})^{\frac{1}{\tilde{\tau}'_j}}(c'_jC_s)^{\frac{1}{\tilde{\tau}_j}}]}
{\sum_{j=1}^2(c_j\|u_0\|^2_{L^2})^{\frac{1}{\tilde{\tau}'_j}}(c'_jC_s)^{\frac{1}{\tilde{\tau}_j}}}
:=\tilde{K}\label{622x1}
\end{align}
and $NF_2(s)s-(N+2)G_2(s)\geq 0$,
then there exists $C$ such that
\begin{align}
\int_{\mathbb{R}^N}[|\nabla h(|u|^2)|^2+|G_1(|u|^2)|+G_2(|u|^2)]dx\leq \frac{C}{t^{2-l_3}}\quad {\rm for}\ t\geq 1,\label{6101'}
\end{align}
where $$l_3=\frac{\tilde{k}_1\sum_{j=1}^2(c_j\|u_0\|^2_{L^2})^{\frac{1}{\tilde{\tau}'_j}}(c'_jC_s)^{\frac{1}{\tilde{\tau}_j}}}
{[1-\sum_{j=1}^2(c_j\|u_0\|^2_{L^2})^{\frac{1}{\tilde{\tau}'_j}}(c'_jC_s)^{\frac{1}{\tilde{\tau}_j}}]}<2.$$

(4) If $2h''(s)h'(s)s+(h'(s))^2\geq 0$, $0<NF_1(s)s-(N+2)G_1(s)\leq\tilde{k}_1|G_1(s)|$ for some $0<\tilde{k}_1<
\tilde{K}$, where $\tilde{K}$ is defined by (\ref{622x1}), and $-k_1G_2(s)\leq NF_2(s)s-(N+2)G_2(s)<0$ for some $0<k_1<2[1-\sum_{j=1}^2(c_j\|u_0\|^2_{L^2})^{\frac{1}{\tilde{\tau}'_j}}(c'_jC_s)^{\frac{1}{\tilde{\tau}_j}}]$,
then there exists $C$ such that
\begin{align}
\int_{\mathbb{R}^N}[|\nabla h(|u|^2)|^2+|G_1(|u|^2)|+G_2(|u|^2)]dx\leq \frac{C}{t^{2-l_4}}\quad {\rm for}\ t\geq 1, \label{6102'}
\end{align}
where
$$
l_4=\frac{\max[k_1, \tilde{k}_1\sum_{j=1}^2(c_j\|u_0\|^2_{L^2})^{\frac{1}{\tilde{\tau}'_j}}(c'_jC_s)^{\frac{1}{\tilde{\tau}_j}}]}
{[1-\sum_{j=1}^2(c_j\|u_0\|^2_{L^2})^{\frac{1}{\tilde{\tau}'_j}}(c'_jC_s)^{\frac{1}{\tilde{\tau}_j}}]}<2.
$$

(5) If $-k_2(h'(s))^2<2h''(s)h'(s)s+(h'(s))^2<0$ for some \begin{align}
0<k_2<\frac{2[1-\sum_{j=1}^2(c_j\|u_0\|^2_{L^2})^{\frac{1}{\tilde{\tau}'_j}}(c'_jC_s)^{\frac{1}{\tilde{\tau}_j}}]}{N}:=K,\label{622x2}
\end{align}
 $NF_1(s)s-(N+2)G_1(s)\leq 0$ and $NF_2(s)s-(N+2)G_2(s)\geq0 $, then there exists $C$ such that
\begin{align}
\int_{\mathbb{R}^N}[|\nabla h(|u|^2)|^2+|G_1(|u|^2)|+G_2(|u|^2)]dx\leq \frac{C}{t^{2-l_5}}\quad {\rm for}\ t\geq 1,\label{6103}
\end{align}
where
$$
l_5=\frac{Nk_2}{[1-\sum_{j=1}^2(c_j\|u_0\|^2_{L^2})^{\frac{1}{\tilde{\tau}'_j}}(c'_jC_s)^{\frac{1}{\tilde{\tau}_j}}]}<2.
$$

(6) If $-k_2(h'(s))^2<2h''(s)h'(s)s+(h'(s))^2<0$ for some $0<k_2\leq K$, where $K$ is defined by (\ref{622x2}), $NF_1(s)s-(N+2)G_1(s)\leq 0$ and $-k_1G_2(s)\leq NF_2(s)s-(N+2)G_2(s)<0$ for some $0<k_1<2[1-\sum_{j=1}^2(c_j\|u_0\|^2_{L^2})^{\frac{1}{\tilde{\tau}'_j}}(c'_jC_s)^{\frac{1}{\tilde{\tau}_j}}]$, then there exists $C$ such that
\begin{align}
\int_{\mathbb{R}^N}[|\nabla h(|u|^2)|^2+|G_1(|u|^2)|+G_2(|u|^2)]dx\leq \frac{C}{t^{2-l_6}}\quad {\rm for}\ t\geq 1,\label{6104}
\end{align}
where
$$
l_6=\frac{\max(Nk_2,k_1)}{[1-\sum_{j=1}^2(c_j\|u_0\|^2_{L^2})^{\frac{1}{\tilde{\tau}'_j}}(c'_jC_s)^{\frac{1}{\tilde{\tau}_j}}]}<2.
$$

(7)  If $-k_2(h'(s))^2<2h''(s)h'(s)s+(h'(s))^2<0$ for some $0<k_2\leq K$, where $K$ is defined by (\ref{622x2}), $0<NF_1(s)s-(N+2)G_1(s)\leq\tilde{k}_1|G_1(s)|$ for some $0<\tilde{k}_1<2$ and $NF_2(s)s-(N+2)G_2(s)\geq 0$,
moreover,
$$
l_7=\frac{[Nk_2+\tilde{k}_1\sum_{j=1}^2(c_j\|u_0\|^2_{L^2})^{\frac{1}{\tilde{\tau}'_j}}(c'_jC_s)^{\frac{1}{\tilde{\tau}_j}}]}
{[1-\sum_{j=1}^2(c_j\|u_0\|^2_{L^2})^{\frac{1}{\tilde{\tau}'_j}}(c'_jC_s)^{\frac{1}{\tilde{\tau}_j}}]},
$$
 then there exists $C$ such that
\begin{align}
\int_{\mathbb{R}^N}[|\nabla h(|u|^2)|^2+|G_1(|u|^2)|+G_2(|u|^2)]dx\leq \frac{C}{t^{2-l_7}}\quad {\rm for}\ t\geq 1.\label{61213}
\end{align}

(8)  If $-k_2(h'(s))^2<2h''(s)h'(s)s+(h'(s))^2<0$ for some $0<k_2\leq  K$, where $K$ is defined by (\ref{622x2}), $0<NF_1(s)s-(N+2)G_1(s)\leq\tilde{k}_1|G_1(s)|$ for some $0<\tilde{k}_1<
\tilde{K}$, where $\tilde{K}$ is defined by (\ref{622x1}), and $-k_1G_2(s)\leq NF_2(s)s-(N+2)G_2(s)<0$ for some $0<k_1<2[1-\sum_{j=1}^2(c_j\|u_0\|^2_{L^2})^{\frac{1}{\tilde{\tau}'_j}}(c'_jC_s)^{\frac{1}{\tilde{\tau}_j}}]$,
moreover,
$$
l_8=\frac{\max[k_1,Nk_2+\tilde{k}_1\sum_{j=1}^2(c_j\|u_0\|^2_{L^2})^{\frac{1}{\tilde{\tau}'_j}}(c'_jC_s)^{\frac{1}{\tilde{\tau}_j}}]}
{[1-\sum_{j=1}^2(c_j\|u_0\|^2_{L^2})^{\frac{1}{\tilde{\tau}'_j}}(c'_jC_s)^{\frac{1}{\tilde{\tau}_j}}]},
$$
then there exists $C$ such that
\begin{align}
\int_{\mathbb{R}^N}[|\nabla h(|u|^2)|^2+|G_1(|u|^2)|+G_2(|u|^2)]dx\leq \frac{C}{t^{2-l_8}}\quad {\rm for}\ t\geq 1. \label{61214}
\end{align}
In all cases above, we have
\begin{align}
\lim_{t\rightarrow \infty}\int_{\mathbb{R}^N}|\nabla u(\cdot,t)|^2dx=2E(u_0),\quad \lim_{t\rightarrow \infty}\|u(\cdot,t)\|^2_{H^1}= M(u_0)+2E(u_0).\label{10281}
\end{align}

2. Assume that $u$ is the blowup solution of (\ref{1}),  $u_0\in X$ and $xu_0\in L^2(\mathbb{R}^N)$. If $[(h'(s))^2+2h''(s)h'(s)s]\leq 0]$, $NF(s)s-(N+2)G(s)\geq 0$ and $-4T^2E(u_0)-\int_{\mathbb{R}^N}|xu_0|^2dx-4T\Im\int_{\mathbb{R}^N}\bar{u}_0(x\cdot \nabla u_0)dx>0$, then
\begin{align}
\int_{\mathbb{R}^N}G(|u|^2)dx\geq \frac{C}{(T-t)^2}.\label{895}
\end{align}
Consequently,
\begin{align}
\int_{\mathbb{R}^N}|\nabla u(\cdot,t)|^2dx+\int_{\mathbb{R}^N}|\nabla h(|u(\cdot,t)|^2)|^2dx\geq \frac{C}{(T-t)^2} \quad {\rm as}\ t \ {\rm close} \ T.\label{10282}
\end{align}

}

\vskip .2in

Before we prove Theorem 5, we would like to recall the following Gronwall's inequality in differential form:

{\bf Gronwall's inequality} {\it Let $\xi(t)$ be a nonnegative, absolutely continuous function on $[a,+\infty)$, which satisfies
\begin{align}
\xi'(t)\leq \phi(t)\xi(t)+\psi(t),\label{624x1}
\end{align}
where $\phi(t)$ and $\psi(t)$ are nonnegative, summable functions on $[a, +\infty)$. Then
\begin{align}
\xi(t)\leq e^{\int_a^t\phi(\eta)d\eta}[\xi(a)+\int_a^t\psi(\eta)e^{-\int_a^{\eta}\phi(\xi)d\xi}d\eta]\label{624x2}
\end{align}
for all $t\in [a, +\infty)$.}

We also point out a fact below
\begin{align}
\int_{\mathbb{R}^N}|G_1(|u|^2)|dx
&\leq \sum_{j=1}^2(c_j\|u_0\|^2_{L^2})^{\frac{1}{\tilde{\tau}'_j}}(c'_jC_s)^{\frac{1}{\tilde{\tau}_j}}\int_{\mathbb{R}^N}|\nabla h(|u|^2)|^2dx,\label{628x1}
\end{align}
which means that
\begin{align}
-4t^2\int_{\mathbb{R}^N}G_1(|u|^2)dx\geq -4t^2\sum_{j=1}^2(c_j\|u_0\|^2_{L^2})^{\frac{1}{\tilde{\tau}'_j}}(c'_jC_s)^{\frac{1}{\tilde{\tau}_j}}\int_{\mathbb{R}^N}|\nabla h(|u|^2)|^2dx.\label{6231}
\end{align}

{\bf Proof of Theorem 5:} 1. Assume that $u$ is the global solution of (\ref{1}),  $u_0\in X$ and $xu_0\in L^2(\mathbb{R}^N)$.

(1) $2h''(s)h'(s)s+(h'(s))^2\geq 0$, $NF_1(s)s-(N+2)G_1(s)\leq 0$ and $NF_2(s)s-(N+2)G_2(s)\geq 0 $.

(\ref{691}), (\ref{628x1}) and (\ref{6231}) imply that
$$
\int_{\mathbb{R}^N}[|\nabla h(|u|^2)|^2+G_2(|u|^2)]dx\leq Ct^{-2}.
$$
Consequently,
$$\int_{\mathbb{R}^N}|G_1(|u|^2)|dx\leq C\int_{\mathbb{R}^N}|\nabla h(|u|^2)|^2dx\leq Ct^{-2},$$
(\ref{6101}) holds.

(2) $2h''(s)h'(s)s+(h'(s))^2\geq 0$, $NF_1(s)s-(N+2)G_1(s)\leq 0$ and $-k_1G_2(s)\leq NF_2(s)s-(N+2)G_2(s)<0$ for some $0<k_1<2$.

By (\ref{691}), (\ref{628x1}) and (\ref{6231}), we have
\begin{align}
&\quad 4[1-\sum_{j=1}^2(c_j\|u_0\|^2_{L^2})^{\frac{1}{\tilde{\tau}'_j}}(c'_jC_s)^{\frac{1}{\tilde{\tau}_j}}]t^2\int_{\mathbb{R}^N}|\nabla h(|u|^2)|^2dx+4t^2\int_{\mathbb{R}^N}G_2(|u|^2)dx\nonumber\\
&\leq\int_{\mathbb{R}^N}|xu_0|^2dx+4k_1\int_0^t\eta [\int_{\mathbb{R}^N}G_2(|u|^2)dx] d\eta.\label{6232}
\end{align}
Let
$$
A_2(t):=4\int_0^t\eta [\int_{\mathbb{R}^N}G_2(|u|^2)dx] d\eta.
$$
(\ref{6232}) implies
$$
A'_2(t)\leq \frac{C_0}{t}+\frac{k_1}{t}A_2(t).
$$
Using Gronwell's inequality, we have
$$
A_2(t)\leq t^{k_1}[A_2(1)+C-\frac{C}{t^{k_1}}]\leq C't^{k_1}.
$$
Using (\ref{628x1}) and (\ref{6231}) again,  we have
\begin{align*}
&\quad 4[1-(c_1\|u_0\|^2_{L^2})^{\frac{1}{\tau'}}(c_2C_s)^{\frac{1}{\tau}}]t^2\int_{\mathbb{R}^N}|\nabla h(|u|^2)|^2dx+4t^2\int_{\mathbb{R}^N}G_2(|u|^2)dx\nonumber\\
&\leq C_0+C't^{k_1}\leq Ct^{k_1}
\end{align*}
and
\begin{align*}
\int_{\mathbb{R}^N}|\nabla h(|u|^2)|^2dx+\int_{\mathbb{R}^N}G_2(|u|^2)dx\leq \frac{C}{t^{2-k_1}},\quad
\int_{\mathbb{R}^N}|G_1(|u|^2)|dx\leq \frac{C}{t^{2-k_1}},
\end{align*}
(\ref{6102}) holds.

(3) $2h''(s)h'(s)s+(h'(s))^2\geq 0$, $0<NF_1(s)s-(N+2)G_1(s)\leq\tilde{k}_1|G_1(s)|$ for some $
0<\tilde{k}_1<\tilde{K}$ and $NF_2(s)s-(N+2)G_2(s)\geq 0$. Similar to (\ref{6232}), letting
$$
A_3(t):=4\int_0^t\eta[\int_{\mathbb{R}^N}|\nabla h(|u|^2)|^2dx]d\eta.
$$
we have
$$
A'_3(t)\leq \frac{C_0}{[1-\sum_{j=1}^2(c_j\|u_0\|^2_{L^2})^{\frac{1}{\tilde{\tau}'_j}}(c'_jC_s)^{\frac{1}{\tilde{\tau}_j}}]t}
+\frac{\sum_{j=1}^2(c_j\|u_0\|^2_{L^2})^{\frac{1}{\tilde{\tau}'_j}}(c'_jC_s)^{\frac{1}{\tilde{\tau}_j}}\tilde{k}_1}
{[1-\sum_{j=1}^2(c_j\|u_0\|^2_{L^2})^{\frac{1}{\tilde{\tau}'_j}}(c'_jC_s)^{\frac{1}{\tilde{\tau}_j}}]t}
A_3(t).
$$
Using Gronwell's inequality, we get
\begin{align*}
A_3(t)&\leq [A_3(1)+C']t^{\left(\frac{\sum_{j=1}^2(c_j\|u_0\|^2_{L^2})^{\frac{1}{\tilde{\tau}'_j}}(c'_jC_s)^{\frac{1}{\tilde{\tau}_j}}\tilde{k}_1}
{[1-\sum_{j=1}^2(c_j\|u_0\|^2_{L^2})^{\frac{1}{\tilde{\tau}'_j}}(c'_jC_s)^{\frac{1}{\tilde{\tau}_j}}]}\right)}
:=Ct^{l_3}
\end{align*}
and
$$
\int_{\mathbb{R}^N}|\nabla h(|u|^2)|^2dx\leq \frac{C}{t^2}+\frac{Ct^{l_3}}{t^2}\leq \frac{C'}{t^{2-l_3}}.
$$
Then we have
$$
\int_{\mathbb{R}^N}G_2(|u|^2)dx\leq \frac{C}{t^2}+\frac{Ct^{l_3}}{t^2}\leq \frac{C'}{t^{2-l_3}},\quad \int_{\mathbb{R}^N}|G_1(|u|^2)|dx\leq  \frac{C'}{t^{2-l_3}},
$$
and (\ref{6101'}) holds.

(4) $2h''(s)h'(s)s+(h'(s))^2\geq 0$, $0<NF_1(s)s-(N+2)G_1(s)\leq\tilde{k}_1|G_1(s)|$ for some $0<\tilde{k}_1<
\tilde{K}$, where $\tilde{K}$ is defined by (\ref{622x1}), and $-k_1G_2(s)\leq NF_2(s)s-(N+2)G_2(s)<0$ for some $0<k_1<2[1-(c_1\|u_0\|^2_{L^2})^{\frac{1}{\tau'}}(c_2C_s)^{\frac{1}{\tau}}]$.

Similar to Case (2) and (3), letting
$$
A_4(t)=4\int_0^t\eta[\int_{\mathbb{R}^N}[|\nabla h(|u|^2)|^2+G_2(|u|^2)]dx]d\eta,
$$
computing $A'_4(t)$ and using Gronwall's inequality, we get
$$
\int_{\mathbb{R}^N}[|\nabla h(|u|^2)|^2+G_2(|u|^2)]dx\leq \frac{C}{t^2}+\frac{Ct^{l_4}}{t^2}\leq \frac{C'}{t^{2-l_4}}
$$
and
$$
\int_{\mathbb{R}^N}|G_1(|u|^2)|dx\leq C\int_{\mathbb{R}^N}|\nabla h(|u|^2)|^2dx\leq \frac{C'}{t^{2-l_4}},
$$
(\ref{6102'}) holds.

(5) $-k_2(h'(s))^2<2h''(s)h'(s)s+(h'(s))^2<0$ for some $0<k_2<K$,
 $NF_1(s)s-(N+2)G_1(s)\leq 0$ and $NF_2(s)s-(N+2)G_2(s)\geq 0 $.

Similar to the cases above,  letting
$$
A_5(t)=4\int_0^t\eta[\int_{\mathbb{R}^N}|\nabla h(|u|^2)|^2dx]d\eta,
$$
and computing $A'_5(t)$, using Gronwall's inequality, we get
$$
\int_{\mathbb{R}^N}|\nabla h(|u|^2)|^2 dx\leq \frac{C}{t^2}+\frac{Ct^{l_5}}{t^2}\leq \frac{C'}{t^{2-l_5}},
$$
and
$$
\int_{\mathbb{R}^N}|G_1(|u|^2)|dx\leq \frac{C}{t^{2-l_5}},\quad \int_{\mathbb{R}^N}G_2(|u|^2)dx\leq \frac{C}{t^{2-l_5}},
$$
and (\ref{6103}) holds.

(6) $-k_2(h'(s))^2<2h''(s)h'(s)s+(h'(s))^2<0$ for some $0<k_2\leq K$, where $K$ is defined by (\ref{622x2}), $NF_1(s)s-(N+2)G_1(s)\leq 0$ and $-k_1G_2(s)\leq NF_2(s)s-(N+2)G_2(s)<0$ for some $0<k_1<2[1-\sum_{j=1}^2(c_j\|u_0\|^2_{L^2})^{\frac{1}{\tilde{\tau}'_j}}(c'_jC_s)^{\frac{1}{\tilde{\tau}_j}}]$.

Similarly, letting
$$
A_6(t)=4\int_0^t\eta[\int_{\mathbb{R}^N}[|\nabla h(|u|^2)|^2+G_2(|u|^2)]dx]d\eta,
$$
and computing $A'_6(t)$, using Gronwall's inequality, we have
$$
\int_{\mathbb{R}^N}[|\nabla h(|u|^2)|^2+G_2(|u|^2)] dx\leq\frac{C}{t^2}+\frac{Ct^{l_6}}{t^2}\leq \frac{C'}{t^{2-l_6}}
$$
and
$$
\int_{\mathbb{R}^N}|G_1(|u|^2)| dx\leq \frac{C'}{t^{2-l_6}},
$$
(\ref{6104}) holds.

(7) $-k_2(h'(s))^2<2h''(s)h'(s)s+(h'(s))^2<0$ for some $0<k_2\leq  K$, where $K$ is defined by (\ref{622x2}), $0<NF_1(s)s-(N+2)G_1(s)\leq\tilde{k}_1|G_1(s)|$ for some $0<\tilde{k}_1<2$ and $NF_2(s)s-(N+2)G_2(s)\geq 0$.

Similarly, letting
$$
A_7(t)=4\int_0^t\eta[\int_{\mathbb{R}^N}|\nabla h(|u|^2)|^2dx]d\eta
$$
and computing $A'_7(t)$, using Gronwell's inequality, we get
$$
\int_{\mathbb{R}^N}|\nabla h(|u|^2)|^2dx\leq \frac{C}{t^2}+ \frac{Ct^{l_7}}{t^2}\leq \frac{C'}{t^{2-l_7}}
$$
and
$$
\int_{\mathbb{R}^N}|G_1(|u|^2)|dx\leq \frac{C}{t^{2-l_7}},\quad \int_{\mathbb{R}^N}G_2(|u|^2)dx\leq \frac{C}{t^{2-l_7}},
$$
(\ref{61213}) holds.

(8) $-k_2(h'(s))^2<2h''(s)h'(s)s+(h'(s))^2<0$ for some $0<k_2\leq  K$, where $K$ is defined by (\ref{622x2}), $0<NF_1(s)s-(N+2)G_1(s)\leq\tilde{k}_1|G_1(s)|$ for some $0<\tilde{k}_1<\tilde{K}$, where $\tilde{K}$ is defined by (\ref{622x1}), and $-k_1G_2(s)\leq NF_2(s)s-(N+2)G_2(s)<0$ for some $0<k_1<2[1-\sum_{j=1}^2(c_j\|u_0\|^2_{L^2})^{\frac{1}{\tilde{\tau}'_j}}(c'_jC_s)^{\frac{1}{\tilde{\tau}_j}}]$.

Using (\ref{691}), (\ref{628x1}) and (\ref{6231}), we have
 \begin{align}
&\quad 4[1-\sum_{j=1}^2(c_j\|u_0\|^2_{L^2})^{\frac{1}{\tilde{\tau}'_j}}(c'_jC_s)^{\frac{1}{\tilde{\tau}_j}}]t^2\int_{\mathbb{R}^N}|\nabla h(|u|^2)|^2dx+4t^2\int_{\mathbb{R}^N}G_2(|u|^2)dx\nonumber\\
&\leq\int_{\mathbb{R}^N}|xu_0|^2dx+4Nk_2\int_0^t\eta \int_{\mathbb{R}^N}|\nabla h(|u|^2)|^2dx d\eta+4\tilde{k}_1\int_0^t\eta \int_{\mathbb{R}^N}|G_1(|u|^2)|dx d\eta\nonumber\\
&\quad+4k_1\int_0^t\eta \int_{\mathbb{R}^N}G_2(|u|^2)dx d\eta\nonumber\\
&\leq C_0+4\max[k_1, Nk_2+\tilde{k}_1\sum_{j=1}^2(c_j\|u_0\|^2_{L^2})^{\frac{1}{\tilde{\tau}'_j}}(c'_jC_s)^{\frac{1}{\tilde{\tau}_j}}]\int_0^t\eta \int_{\mathbb{R}^N}[|\nabla h(|u|^2)|^2+G_2(|u|^2)]dx.\label{624x3}
\end{align}
Similarly, letting
$$
A_8(t)=4\int_0^t\eta[\int_{\mathbb{R}^N}[|\nabla h(|u|^2)|^2+G_2(|u|^2)]dx]d\eta
$$
and computing $A'_8(t)$, using Gronwell's inequality, we get
$$
\int_{\mathbb{R}^N}[|\nabla h(|u|^2)|^2+G_2(|u|^2)]dx\leq \frac{C}{t^2}+\frac{Ct^{l_8}}{t^2}\leq \frac{C'}{t^{2-l_8}}
$$
and
$$
\int_{\mathbb{R}^N}|G_1(|u|^2)|dx\leq \frac{C}{t^{2-l_8}},
$$
(\ref{61214}) holds.

In all cases above, we have
$$
\lim_{t\rightarrow \infty}\int_{\mathbb{R}^N}|\nabla h(|u|^2)|^2dx=0,\quad \lim_{t\rightarrow \infty}\int_{\mathbb{R}^N}|G_1(|u|^2)|dx=0,\quad
\lim_{t\rightarrow \infty}\int_{\mathbb{R}^N}G_2(|u|^2)dx=0.
$$
Using energy conservation law $E(u)=E(u_0)$ and letting $t\rightarrow \infty$, we get
$$
\frac{1}{2}\lim_{t\rightarrow \infty}\int_{\mathbb{R}^N}|\nabla u|^2dx=\lim_{t\rightarrow \infty}E(u)=E(u_0).
$$
Hence
$$
\lim_{t\rightarrow \infty}\|u\|_{H^1}^2=\lim_{t\rightarrow \infty}\left(\int_{\mathbb{R}^N}|u|^2dx+\int_{\mathbb{R}^N}|\nabla u|^2dx\right)=M(u_0)+2E(u_0),
$$
which prove (\ref{10281}).

2. Assume that $u$ is the blowup solution of (\ref{1}),  $u_0\in X$ and $xu_0\in L^2(\mathbb{R}^N)$.
Using (\ref{893}), we have
\begin{align}
&\quad 4(T-t)^2\int_{\mathbb{R}^N}G(|u|^2)dx\nonumber\\
&=\int_{\mathbb{R}^N}|(x+2i(T-t)\nabla)u|^2dx+4(T-t)^2\int_{\mathbb{R}^N}|\nabla h(|u|^2)|^2dx\nonumber\\
&\quad+4\int_0^t(T-\tau)\theta(\tau)d\tau -4T^2E(u_0)-\int_{\mathbb{R}^N}|xu_0|^2dx\nonumber\\
&\quad-4T\Im\int_{\mathbb{R}^N}\bar{u}_0(x\cdot \nabla u_0)dx.\label{894}
\end{align}
If $[(h'(s))^2+2h''(s)h'(s)s]\leq 0]$, $NF(s)s-(N+2)G(s)\geq 0$ and $$E(u_0)\leq 0,\quad -4T^2E(u_0)-\int_{\mathbb{R}^N}|xu_0|^2dx-4T\Im\int_{\mathbb{R}^N}\bar{u}_0(x\cdot \nabla u_0)dx>0,$$ then (\ref{894}) implies that
\begin{align*}
\int_{\mathbb{R}^N}G(|u|^2)dx\geq \frac{C}{(T-t)^2}.
\end{align*}
(\ref{895}) holds.

Using energy conservation law $E(u)=E(u_0)$, we have
$$
\frac{1}{2}\int_{\mathbb{R}^N}|\nabla u|^2dx+\frac{1}{2}\int_{\mathbb{R}^N}|\nabla h(|u|^2)|^2dx=\frac{1}{2}\int_{\mathbb{R}^N}G(|u|^2)dx+E(u_0)\geq
\frac{C}{(T-t)^2}+E(u_0).
$$
As $t$ close to $T$ enough, then there exists a constant $0<C'<C$ such that
$$
\frac{C}{(T-t)^2}+E(u_0)\geq \frac{C'}{(T-t)^2}.
$$
Consequently,
$$
\int_{\mathbb{R}^N}|\nabla u|^2dx+\int_{\mathbb{R}^N}|\nabla h(|u|^2)|^2dx\geq \frac{2C'}{(T-t)^2},
$$
which proves (\ref{10282}). \hfill $\Box$

\end{document}